\newif\ifpdf
\ifx\pdfoutput\undefined
\pdffalse 
\else
\pdfoutput=1 
\pdftrue \fi

\newif\iffinal
\finaltrue 

\documentclass[reqno,twoside,11pt]{amsart}
\iffinal\else\usepackage[notref,notcite]{showkeys}\fi
\usepackage{amsmath}
\usepackage{amsfonts}
\usepackage{amssymb}
\usepackage{verbatim}
\usepackage{epsfig}
\usepackage{setspace}

\IfFileExists{epsf.def}{\input epsf.def}{\usepackage{epsf}}

\IfFileExists{myowntimes.sty}{\usepackage{myowntimes}\usepackage{mathrsfs}}
	{\usepackage{mathpazo}\usepackage{mathrsfs}}

\DeclareFontFamily{OT1}{eusb}{} \DeclareFontShape{OT1}{eusb}{m}{n} {<5> <6> <7> <8> <9> <10> <11> <12> <14.4> eusb10}{}
\DeclareMathAlphabet{\eusb}{OT1}{eusb}{m}{n}

\DeclareFontFamily{OT1}{eusm}{} \DeclareFontShape{OT1}{eusm}{m}{n} {<5> <6> <7> <8> <9> <10> <11> <12> <14.4> eusm10}{}
\DeclareMathAlphabet{\eusm}{OT1}{eusm}{m}{n}

\DeclareFontFamily{OT1}{eufm}{} \DeclareFontShape{OT1}{eufm}{m}{n} {<5> <6> <7> <8> <9> <10> <11> <12> <14.4> eufm10}{}
\DeclareMathAlphabet{\mathfrak}{OT1}{eufm}{m}{n}

\DeclareFontFamily{OT1}{fraktura}{}
\DeclareFontShape{OT1}{fraktura}{m}{n} {<5> <6> <7> <8> <9> <10> <11> <12> <13> <14.4> [1.1] eufm10}{}
\DeclareMathAlphabet{\fraktura}{OT1}{fraktura}{m}{n}

\DeclareFontFamily{OT1}{cmfi}{} \DeclareFontShape{OT1}{cmfi}{m}{n} {<5> <6> <7> <8> <9> <10> <11> <12> <13> <14.4> [0.9] cmfi10}{}
\DeclareMathAlphabet{\cmfi}{OT1}{cmfi}{b}{n}

\DeclareFontFamily{OT1}{cmss}{} \DeclareFontShape{OT1}{cmss}{m}{n} {<5> <6> <7> <8> <9> <10> <11> <12> <13> <14.4> cmss10}{}
\DeclareMathAlphabet{\cmss}{OT1}{cmss}{m}{n}

\newtheoremstyle{thm}{1.5ex}{1.5ex}{\itshape\rmfamily}{} {\bfseries\rmfamily}{}{2ex}{}

\newtheoremstyle{def}{1.5ex}{1.5ex}{\slshape\rmfamily}{} {\bfseries\rmfamily}{}{2ex}{}

\newtheoremstyle{rem}{1.3ex}{1.3ex}{\rmfamily}{} {\itshape}
{} {1.5ex}{}

\theoremstyle{thm}
\newtheorem{theorem}{Theorem}[section]
\newtheorem{lemma}[theorem]{Lemma}

\newtheorem*{Main Theorem}{Main Theorem.}
\newtheorem{corollary}[theorem]{Corollary}
\newtheorem*{special theorem}{Lindeberg-Feller Theorem for Martingales}

\theoremstyle{def}

\theoremstyle{rem}
\newtheorem{remark}{{\itshape Remark}}[]

\numberwithin{equation}{section}

\renewcommand{\section}{\secdef\sct\sect}
\newcommand{\sct}[2][default]{%
\refstepcounter{section}
\addcontentsline{toc}{section}{{\tocsection {}{\thesection}{\!\!\!\!#1\dotfill}}{}}
\vspace{0.7cm}
\centerline{\scshape\thesection.\ #1} \nopagebreak \vspace{0.2cm}}
\newcommand{\sect}[1]{%
\vspace{0.4cm} \centerline{\large\scshape\rmfamily #1}
\vspace{0.2cm}}

\renewcommand{\subsection}{\secdef\subsct\sbsect}
\newcommand{\subsct}[2][default]{\refstepcounter{subsection}
\addcontentsline{toc}{subsection}
{{\tocsection{\!\!}{\hspace{1.2em}\thesubsection}{\!\!\!\!#1\dotfill}}{}}
\nopagebreak\vspace{0.45\baselineskip} {\flushleft\bf
\thesubsection~\bf #1.~}
\\*[3mm]\noindent
\nopagebreak}
\newcommand{\sbsect}[1]{\vspace{0.1cm}\noindent
\textbf{#1.~}\vspace{0.1cm}}

\renewcommand{\subsubsection}{%
\secdef \subsubsect\sbsbsect}
\newcommand{\subsubsect}[2][default]{%
\refstepcounter{subsubsection} 
\addcontentsline{toc}{subsubsection}{{\tocsection{\!\!}
{\hspace{3.05em}\thesubsubsection}{\!\!\!\!#1\dotfill}}{}}
\nopagebreak
\vspace{0.15\baselineskip} \nopagebreak {\flushleft\rmfamily
\itshape\thesubsubsection
\ \rmfamily #1\/.}\ }
\newcommand{\sbsbsect}[1]{\vspace{0.1cm}\noindent
\rmfamily \itshape
\arabic{section}.\arabic{subsection}.\arabic{subsubsection} \
\sffamily #1\/.\ }

\iffinal

\else

\fi

\renewcommand{\caption}[1]{%
\vglue0.5cm
\refstepcounter{figure}
\begin{minipage}{0.9\textwidth}\small {\sc Figure~\thefigure. }#1\end{minipage}}

\newcommand{\E}{\mathbb E}

\newcommand{\N}{\mathbb N}

\def\myffrac#1#2 in #3{\raise 2.6pt\hbox{$#3 #1$}\mkern-1.5mu\raise 0.8pt\hbox{$#3/$}\mkern-1.1mu\lower 1.5pt\hbox{$#3 #2$}}
\newcommand{\ffrac}[2]{\mathchoice%
	{\myffrac{#1}{#2} in \scriptstyle}
	{\myffrac{#1}{#2} in \scriptstyle}
	{\myffrac{#1}{#2} in \scriptscriptstyle}
	{\myffrac{#1}{#2} in \scriptscriptstyle}
}

\title[CLT for Multi-overlaps]{The Interaction between Multi-overlaps in the High Temperature Phase of the Sherrington-Kirkpatrick Spin Glass}
\author[N.~Crawford]{Nicholas Crawford}

\begin{document}
\thanks{\hglue-4.5mm\fontsize{9.6}{9.6}\selectfont\copyright\,2007 by N.~Crawford. Reproduction, by any means, of the entire
article for non-commercial purposes is permitted without charge.\vspace{2mm}}
\maketitle

\vspace{-5mm}
\centerline{\textit{Department of Statistics, University of California at Berkeley}}

\vspace{-2mm}
\begin{abstract}
In this paper we explore the joint behaviour of a finite number of multi-overlaps in the high temperature phase of the SK model.  Extending work by M. Talagrand, we show that, when these objects are scaled to have non-trivial limiting distributions, the joint behaviour is described by a Gaussian process with an explicit covariance structure. 
\end{abstract}

\section{Introduction}
\noindent
The recent text \cite{Talagrand-book} provides a beautiful introduction to interpolation and the cavity method and their use within the realm of mean field spin glasses.  In the particular case of the SK model it provides detailed information about quantities of physical interest.  It was used at high temperatures to compute the quenched free energy of the system precisely \cite{Talagrand-paper-1}.  Moreover the ideas are interwoven with more sophisticated interpolations to obtain upper bounds on the free energy at all temperatures, \cite{Guerra}, and, ultimately, to express  this quantity in terms of a variational principle predicted by G. Parisi \cite{Talagrand-paper}.  

Recall that the SK model of spin glasses is defined as the Gibbs measure on spin configurations $\sigma \in  \{-1, 1\}^N$ with the Hamiltonian
\begin{equation}
H_N(\sigma) = \sum_{1 \leq i<j \leq N} \frac{1}{\sqrt N} g_{i,j} \sigma_i \sigma_j + h \sum_{i=1}^N \sigma_i.
\end{equation}
Here the couplings $g_{i,j}$ are taken to be independent Gaussian random variables with mean $0$ and variance $1$ and $h \in \mathbb R$ is the strength of the external field which may, without loss of generality, be assumed to be positive.  In other words, each spin configuration is chosen with probability 
\begin{equation}
P(\sigma) \propto
e^{\beta H_N(\sigma)}
\end{equation}
where we have omitted the minus sign from the exponent for convenience and the parameter $\beta$ denotes the inverse temperature.

Let $\sigma^1, \boldmath\sigma^2 \in \{-1,1\}^N$ denote a pair of spin configurations and let the overlap between $\sigma^1, \boldmath\sigma^2$ be defined by
\begin{equation}
R_{1,2}= \frac{1}{N}\sum_{i=1}^N \sigma_i^1 \sigma_i^2.
\end{equation}
M. Talagrand \cite{Talagrand-book} showed that for high enough temperature, there exists a value $q_2 \in \mathbb R$ so that 
\begin{equation}
\label{Eq:High-T}
\E \left[\left\langle \left(R_{1,2} - q_2\right)^2 \right\rangle\right]  \leq \frac{L}{N}
\end{equation}
for some constant $L>0$ as $N$, the number of spins in the system, tends to $\infty$.  Here $\langle \cdot \rangle$ is the quenched Gibbs state on $N$ spins and $\nu\left( \cdot \right)\E\left[ \cdot \right]$ denotes the average over the disorder of the system.  

Further, considering the quenched Gibbs state for $n$ replicas, in Sections $2.5$, $2.6$, and $2.7$ of \cite{Talagrand-book}, Talagrand develops the machinery needed to compute the joint distribution of the scaled random variables
\begin{equation}
\{\sqrt N \left(R_{\ell, \ell'}-q_2\right)\}_{\{\ell, \ell'\} \subset \{1, \dotsc, n\}}.
\end{equation}
Here the replicas $\{\sigma_{\ell}\}$ are spin configurations sampled independently according to the quenched Gibbs state with the \textit{same} realization of disorder.
 
Once the control \eqref{Eq:High-T} of the overlaps is obtained, it may be shown via the cavity method that in this high temperature regime, many quantities of interest are concentrated in the above sense as well.  For example, it is of physical interest to ask about the behaviour of the magnetization
\begin{equation}
M= \frac{1}{N} \sum_{i=1}^N \sigma_i
\end{equation}
under the quenched Gibbs state.
This random variable is concentrated, and if we denote by $q_1$ the value that $M$ concentrates around then we can immediately generalize the above and ask about the joint distribution of
\begin{equation}
\{\sqrt N\left(R_{\ell, \ell'}-q_2\right)\}_{\{\ell, \ell'\} \subset \{1, \dotsc, n\}} \quad 
\{\sqrt N\left(M_{\ell}-q_1\right)\}_{\ell \in \{1, \dotsc, n\}}.
\end{equation}

Towards the end of Section 2.7 of \cite{Talagrand-book}, a generalization of this problem is proposed:  Given a subset $S \subset \{1, \dotsc, n\}$ of replica indices, we may introduce  the multi-overlaps
\begin{equation}
R_S = \frac{1}{N} \sum_{i=1}^N \prod_{\ell \in S} \sigma^\ell_i.
\end{equation}
Again, given \eqref{Eq:High-T}, it is not difficult to show that these objects are concentrated around certain values $q_{|S|}$, where $|S|$ denotes the cardinality of $S$.  The general challenge is to characterize the joint distribution of the scaled variables
\begin{equation}
\{\sqrt{N}\left(R_S -q_{|S|}\right)\}_{S \subseteq \{1, \dotsc, n\}}.
\end{equation}
This is the question addressed in the present paper.
We shall not give a precise description here, as it requires the introduction of a bit of notation, however let us simply say that the joint distribution is given by a Gaussian process with an explicit covariance structure.  

Unfortunately there is no magic bullet here (besides the miracles inherent when working with Ising variables).   Our results are proved by the method of moments, identifying the underlying covariance structure first and then inductively computing joint moments.

In the next section, we shall introduce this notation and give the main results in a precise way.  Subsequent sections are devoted to proving these results.
We will proceed with proofs in a general pattern.  It turns out that the characterization of distributions \textit{decrease} in difficulty as $|S|$ increases.  Thus we always proceed from easiest to hardest case (the magnetizations being basically the most difficult).  This creates redundancy, but probably increases readability.

Finally, let us mention some overlapping work related to this result which we became aware of during the preparation of this paper.  Besides the standard on the subject \cite{Talagrand-book} to which we refer below, the paper \cite{Tindel} treats this problem with no external field, however this represents a significant simplification over the general case.  Also, that paper is interested in a detailed expansion of the corrections to computed moments.
Another work which shares some features in common with the present work is \cite{Hanen}, which considers the behavior of the random variables $\langle \sigma_i \sigma_j\rangle$.

\section{Notation and Results}
Throughout the following, the constant $C>0$ will be used to denote a quantity that does not depend on $N$, though it may depend on $\beta$ and the number of replicas in the quantity of interest.  This value will change from instance to instance.  We shall further (following Talagrand) use the notation $O(k)$ to denote a quantity bounded by $\frac{C}{N^{-\frac{k}{2}}}$ in absolute value.

Let us denote the truncated spins by 
\begin{equation}
\dot{\sigma_i}^\ell= \sigma_i^\ell- \langle \sigma_i \rangle
\end{equation} 
where the superscript indicates the replica of interest and the subscript indicates the site of interest.

Let $q_2$ denote the solution to the equation 
\begin{equation}
\label{Eq:q_2}
q_2= \mathbb E\left[ \tanh^2 \left(\beta \sqrt{q_2} Y +h \right)\right]
\end{equation}
where $Y$ is a standard Gaussian.
A result of Guerra and Lata{\l}a \cite{Guerra-2, Latala} shows that there is a unique solution for $q_2$ whenever $h>0$.
Let us further denote 
\begin{equation}
q_p= \mathbb E \left[ \tanh^p \left(\beta \sqrt{q_2} Y +h \right)\right].
\end{equation}

To prove the multi-overlap CLT we need high temperature conditions based on the work of Talagrand.  We assume in what follows that the quantity $N(R_{1,2}-q_2)^2$ has an exponential moment.  This is the content of the following theorem.
\begin{theorem}[Theorem 2.5.1 of \cite{Talagrand-book}]
\label{T:Exp-mom}
There exists $\beta_0>0$ (independent of $h$) such that if $\beta < \beta_0$ we have
\begin{equation}
\nu\left(\exp\left(\frac NL (R_{1,2}-q_2)^2\right)\right) \leq L
\end{equation}
for some $L>0$.
\end{theorem}

For each finite subset $S \subset \N$ and $p \in \N\cup \{0\}$,
let us introduce the following quantity which generalizes truncated overlaps defined by Talagrand in \cite{Talagrand-book}:
\begin{equation}
T_{S,p}= \frac{1}{N} \sum_{i=1}^N  \prod_{\ell \in S} \dot{\sigma_i}^\ell\langle \sigma_i \rangle^p \quad  \text{ when $S \neq \varnothing$}
\end{equation}
and  
\begin{equation}
T_{\varnothing, p}= \frac{1}{N} \sum_{i=1}^N \langle \sigma_i \rangle^p - q_p \quad \text{ when $S = \varnothing$.}
\end{equation}
In the case $p=0$, we shall often denote $T_{S, 0}= T_{S}$.  We remark that this last notation matches Talagrand's definition when $|S|=2$ but not when $|S| \in\{0,1\} $.  Also whe $S=\{\ell\}$, we shall denote
$T_{S, p}= T_{\ell, p}$.  It seems an appropriate time to mention that below the reader will encounter ranges of indices that may lead to empty sums or products.  It should be clear from the context, but we mention here that in all cases empty sums are interpreted as $0$ and empty products are interpreted as $1$.

Our interest in these quantities stems from the fact that we may express the multioverlaps of the SK model in terms of them:
\begin{equation}
\frac{1}{N} \sum_{i=1}^N \sigma_i^1 \cdots \sigma_i^r - q_r = \sum_{S \subseteq [r]} T_{S, r-|S|}.
\end{equation}
Here the notation $[r]$ is used to denote the subset $\{1, \dots, r\}\subset \N$.

To prove the CLT for multi-overlaps, we will first characterize the joint distribution of the quantities
\begin{equation}
\{\sqrt{N} T_{S, p}\}_{S \subset [n], p \in \mathbb N\cup\{0\}}.
\end{equation} 
This characterization goes by the method of moments, facilitated by the fact that whnever $S \vartriangle \tilde S \neq \varnothing$, the families 
\begin{equation}
\{\sqrt{N} T_{S, p}\}_{p \in \mathbb N\cup \{0\}}, \quad \{\sqrt{N} T_{\tilde S, p}\}_{p \in \mathbb N\cup\{0\}} .
\end{equation} are asymptotically finitely independent.

It is notationally convenient for us to introduce a number of quantities defined in terms of the moments $\{q_{p}\}_{p \in \N}$ which come up naturally in the calculations below.
\begin{align}
&\varrho_{p}=q_{p} - 2q_{p+2} + q_{p+4}\\
&\pi_{p}=
\begin{cases}
q_{p+1}-q_{p+3} \hspace{12pt} \text{ if $p \geq -1$}\\
0 \hspace{24pt} \text{ else}
\end{cases}
\\
&\varphi_{p}= p\pi_{p-2} -3\pi_p
\end{align}
Let us further introduce the quantities
\begin{equation}
 A_s\left(p, \tilde p\right) = 
 \begin{cases}
 &\frac{1}{N}\sum_{\ell = 0}^{s} {s \choose \ell} \left(-1\right)^{\ell} q_{p+\tilde{p} + 2 \ell} \quad \quad \quad \quad \text{if $s\geq 3$}.\\
 \\
 &\frac{1}{N} \varrho_{p+\tilde p}+ \frac{\beta^2}{N}\varrho_{p}\varrho_{\tilde p} + \frac{\beta^4}{N}\varrho_{p} \varrho_{\tilde p}  \frac{\varrho_{0}}{1- \beta^2\varrho_{0}} \quad \text{ if $s=2$}.\\
 \\
 &\frac{1}{N}\pi_{p+\tilde{p}-1}+
\beta^2\varphi_p A_1(\tilde p, 1)\\
& + \beta^2\tilde p\pi_{p}A_2(\tilde p - 1, 0) \quad \quad \quad \quad \quad \text{ if $s=1$}\\
 \\
 & \beta^2 \lambda _{p} A_0\left(2, \tilde p\right) 
+\beta^2\kappa_{p, \tilde p} A_1\left(1, \tilde p-1\right) \\
&+\beta^2 \omega_{p, \tilde p}A_2(2, \tilde p-2)
 \quad \quad \quad \text{if $s=0$}
 \end{cases}
 \end{equation}
where
\begin{align*}
\lambda_{p}=& {p \choose 2} \left(q_{p-2}- q_pq_2\right)+ { p+1 \choose 2} \left(q_{p+2}-q_p q_2\right) \\ &- p^2 \left(q_p-q_pq_2\right).\\
\kappa_{p, \tilde p} =&  p\tilde p \left(q_p- q_{p+2}\right) - \tilde p\left(p+1\right) \left(q_{p+2}-q_pq_2\right)\\
\omega_{p, \tilde p}= & {\tilde{p}  \choose 2} \left(q_{p+2}-q_pq_2\right).
\end{align*}
A comment on the cases $s=1,0$ is in order.  Setting $\tilde p =1$, (resp. $\tilde p=2$), we may solve the equation for all $p$ and then apply the result to the solve general case.  
 
The actual values of $A_s(p, \tilde p)$ are less important than their significance: 
\begin{theorem}[The Covariance Structure of the Truncated Family]
\label{T:Cov}
Suppose that $S, \tilde S \subset \N$ are finite subsets.  Then for any fixed $p, \tilde{p} \in \N$,
\begin{equation}
\nu\left(T_{S,p}T_{\tilde{S}, \tilde{p}}\right) = \frac{1}{N}\delta_{S,\tilde S} A_{|S|}\left( p, \tilde p \right) 
\end{equation}
\end{theorem}

\begin{remark}
In the cases $|S|=2$, $p=\tilde{p} = 0$ and $|S|=1$, $p=\tilde p=1$, an easy calculation shows that this formula agrees with the moments computed by Talagrand (as it should...see Theorem \ref{T:CLT-Tal} below). 
\end{remark}

In what follows let $\{Y_{S,p}\}$ be a family of Gaussian random variables independent from the randomness in the SK model with joint covariance structure given by Theorem \ref{T:Cov}.  
As mentioned above, to prove the assertions of this paper, we shall rely on computations due to Talagrand which can be found in Sections 2.5, 2.6 and 2.7 of \cite{Talagrand-book}.  The upshot of these computations is the following.
\begin{theorem}[Theorem 2.7.3 from \cite{Talagrand-book}]
\label{T:CLT-Tal}
Let $n\in \N$ be fixed.   Consider the collection of subsets $S \subseteq [n]$ with $|S| \leq 2$.  Let us fix a collection of integers $k\left(S\right)$ associated to these subsets.  Let $k= \sum_{S\subset [n]:|S| \leq 2} k(S)$.  Then
\begin{equation}
\left|\nu\left(\prod_{S \subseteq [n], \: |S| \leq 2} T_{S}^{k\left(S\right)}\right)- N^{-\ffrac{k}{2}} \prod_{S \subseteq [n]: |S| \leq 2} \E Y_{S}^{k\left(S\right)}\right| \leq \frac{K\left(k\right)}{N^{\frac{k+1}{2}}} 
\end{equation}
\end{theorem}

In light of this result, the main theorem of our work reads as follows.

\begin{theorem}[The Full Picture]
\label{T:Main}
Let $m, n \in \N$ be fixed.
Consider the family $\{T_{S,p}\}_{S\subseteq [n],\:0 \leq  p \leq m}$ and an arbitrary collection of integers $\{k\left(S,p\right)\}_{S\subseteq [n], \:0 \leq  p \leq m}$.  Let $k= \sum_{(S,p)} k\left(S,p\right)$.  Then we have
\begin{equation}
\left|\nu\left(\prod_{(S,p)} T_{S,p}^{k\left(S,p\right)}\right) - N^{-\ffrac{k}{2}} \prod_{S \subseteq [n]} \E \left[\prod_{p=0}^m Y_{S,p}^{k\left(S,p\right)}\right]\right| \leq \frac{C}{N^{\frac{k+1}{2}}} 
\end{equation}
\end{theorem}

There is an  obvious corollary to the previous theorem:
\begin{corollary}
Let $n \in \N$ be fixed.
Then the family $\{\sqrt{N}\left(R_{S}- q_{|S|}\right)\}_{S\subseteq [n]}$ is jointly normal with an explicit covariance structure that can be recovered from Theorem \ref{T:Main} combined with the fact that 
\begin{equation}
R_{S}- q_{|S|}= \sum_{S' \subseteq S} T_{S', |S|-|S'|}
\end{equation}

\end{corollary}

\section{Preliminaries}

The single most important tool in this work appears in \cite{Talagrand-book}.  It is a consequence of the cavity interpolation scheme.  Let us recall briefly the idea.  Let $\Sigma_N = \{-1, 1\}^N$ and let $\Sigma_N^n$ be the space of $n$ replicas.  On $\Sigma_N$ we isolate the final spin and denote it by $\varepsilon$ while in the space of $n$ replicas we denote the $n$ final spins by $\varepsilon^1, \dotsc \varepsilon^n$. 

In the case $n=1$, this final spin feels an effective field from the previous $N-1$ spins which takes the form
\begin{equation}
\rho_N(\sigma)= \frac{1}{\sqrt N}\sum_{i<N} g_{i,N}\sigma_i + h.
\end{equation}
In order to study the effect of this field we introduce the interpolation
\begin{equation}
H_N(t) = H_N^- + \sqrt{t} \rho_N(\sigma)\varepsilon  + \sqrt{1-t} \left(\sqrt{q_2} Y + h\right)\varepsilon
\end{equation}
where $Y$ is a standard Gaussian and $q_2$ solves the equation \eqref{Eq:q_2}.  Let us define $\langle \cdot \rangle_t$ denote the Gibbs state corresponding to this Hamiltonian at inverse temperature $\beta$.  We use the same notation for the Gibbs state on $\Sigma_N^n$ defined by drawing $n$ independent copies of $\sigma$ with respect to the same realization of disorder.  It will be clear from the context which state we are interested in.  Finally we let
\begin{equation}
\nu_t(f)= \E\left[\langle f \rangle_t \right]
\end{equation}
for each $f: \Sigma_N^n \rightarrow \mathbb R$.

The key lemma, proved in \cite{Talagrand-book}, is as follows.
\begin{lemma}[Proposition 2.5.3 from \cite{Talagrand-book}]
\label{L:Cav-Bound}
Suppose $\beta_0$ is as in Theorem \ref{T:Exp-mom}.
Let $ \beta \leq \beta_0$.
Then for any $f$ on $\Sigma^n_N$,
\begin{align}
&\left|\nu\left(f\right)- \nu_{0}\left(f\right)\right| \leq \frac{K\left(n\right)}{\sqrt{N}} \nu\left(f^2\right)^{\frac 12}\\
&\left|\nu\left(f\right)- \nu_{0}\left(f\right) - \nu_0'\left(f\right)\right| \leq \frac{K\left(n\right)}{N} \nu\left(f^2\right)^{\frac 12}.
\end{align}
\end{lemma}

Our main function of interest below will be of the form 
\begin{equation}
\nu\left(T_{S_0,p_0} \mathcal G\right)
\end{equation}
where $\mathcal G$ is a monomial in the variables $\{T_{S,p}\}_{S \subseteq[n], 0 \leq p \leq m}$ and $S_0 \subseteq [n], p_0 \in \N$.  We will need to expand this product in terms of replicas.  We will always use the following notation when making such an expansion.  

To represent $\mathcal G$ in terms of a disjoint set of replicas, consider a total ordering of the factors $T_{S,p}$ appearing in $\mathcal G$.  Denote this ordering by $\mathcal A$.  Recursively, suppose we have introduce replicas up to the index $a\in \mathcal A$.  If $T_{S,p}$ corresponds to the index $a$, define a pair of injections $\alpha_a: S \rightarrow \N$ and $\gamma_a: [p] \rightarrow \N$ so that their images are disjoint from $[n]$ and the images of the previously defined maps.   Finally, we distinguish a pair of injections $\zeta: S_0 \rightarrow \N$ and $\eta: \{0, \dotsc, p\} \rightarrow \N$ whose images are disjoint from each other, all, previous images and $[n]$.  Thus, to be explicit, we may expand the factor $T_{S,p}$ as
\begin{equation}
T_{S,p}= \frac{1}{N}\sum_{i=1}^N \prod_{r \in S}\left(\sigma_i^{r}-\sigma_i^{\alpha_a(r)}\right) \sigma_i^{\gamma_a(1)} \cdots \sigma_i^{\gamma_a(p)}
\end{equation}
Note that for each replica $\ell \notin [n]$, there is a \textit{unique} factor $T_{S,p}$ of $\mathcal G$ which depends on this replica.  We shall denote this dependence by $T_{S,p}(\ell)$.

We record a general computation that we shall use repeatedly in various contexts below:
\begin{lemma}
\label{L:cavitation}
Let $S,p$ be fixed.  Let $n$ be large enough so that $S \cup \zeta\left(S\right)\cup \eta \left([p]\right) \subseteq [n]$.  Suppose $f^{-}: \Sigma_{N-1}^n \rightarrow \mathbb R$.  Then
\begin{equation}
\begin{split}
&\nu\left(\prod_{r\in S} \left(\varepsilon^{r}- \varepsilon^{\zeta\left(r\right)}\right) \prod_{l=1}^p \varepsilon^{\eta\left(l\right)}f^{-}\right)= O\left(2\right)\nu\left(\left(f^{-}\right)^2\right)^{\frac 12}  \hspace{32pt} \text{ if $|S| \geq 3$}.\\
&\hspace{16pt} = \beta^2\varrho_p \nu_{0}\left(f^-\left(R_{r_1, r_2}^--R_{r_1, \zeta\left(r_2\right)}^- - R_{\zeta\left(r_1\right), r_2}^-+R_{\zeta\left(r_1\right), \zeta\left(r_2\right)}^-\right)\right)\\
&  \hspace{48pt} + O\left(2\right)\nu\left(\left(f^{-}\right)^2\right)^{\frac 12} \hspace{130pt} \text{if $S=\{r_1, r_2\}$}.\\
& \hspace{16pt}= \beta^2 \pi_{p-2}\sum_{1 \leq l \leq p} \nu_{0}\left(f^-\left(R_{r_1,\eta\left(l\right)}^--R_{\zeta\left(r_1\right), \eta\left(l\right)}^-\right)\right)\\
&\hspace{48pt} + \beta^2 \pi_{p}\sum_{ \ell \notin\{r_1,\zeta\left(r_1\right)\}\cup \eta\left([p]\right)} \nu_{0}\left(f^-\left(R_{r_1, \ell}^--R_{\zeta\left(r_1\right), \ell}^-\right)\right)\\
& \hspace{48pt} - \beta^2 n \pi_{p}\nu_{0}\left(f^-\left(R_{r_1, n+1}^--R_{\zeta\left(r_1\right), n+1}^-\right)\right)\\
&\hspace{48pt} + O\left(2\right)\nu\left(\left(f^{-}\right)^2\right)^{\frac 12}  \hspace{130pt} \text{if $S = \{r_1\}$}.
\end{split}
\end{equation}

Finally, if $|S|=0$,
\begin{multline}
\nu\left(\left(\prod_{l=1}^p \varepsilon^{\eta\left(l\right)}- q_p\right)f^{-}\right)=
\beta^2 \left(q_{p-2}- q_pq_2\right)\sum_{\{\ell, \ell'\} \subseteq [p]} \nu_{0}\left(f^-\left(R_{\ell, \ell'}^--q_2\right)\right)\\
+\beta^2 \left(q_{p+2}- q_p q_2\right)\sum_{\{\ell, \ell'\}  \cap  \eta([p]) = \varnothing} \nu_{0}\left(f^-\left(R_{\ell, \ell'}^--q_2\right)\right)\\
+\beta^2 \left(q_p- q_p q_2\right)\sum_{\stackrel{1 \leq \ell \leq p}{p< \ell' \leq n}} \nu_{0}\left(f^-\left(R_{\ell, \ell'}^--q_2\right)\right)\\
- \beta^2 n \left(q_p-q_pq_2\right)\sum_{1 \leq \ell \leq p} \nu_{0}\left(f^-\left(R_{\ell, n+1}^--q_2\right)\right)\\
- \beta^2 n \left(q_{p+2}-q_pq_2\right)\sum_{p< \ell \leq n} \nu_{0}\left(f^-\left(R_{\ell, n+1}^--q_2\right)\right)\\
+  \beta^2 \frac{n\left(n+1\right)}{2} \left(q_{p+2}-q_p q_2\right)\nu_{0}\left(f^-\left(R_{n+1, n+2}^--q_2\right)\right) + O(2) \nu\left((f^-)^2\right)^{\frac{1}{2}}.
\end{multline}
\end{lemma}

\begin{proof}
Lemma \ref{L:Cav-Bound} asserts that for any $f$ on $\Sigma^n_N$,
\begin{equation}
\left|\nu\left(f\right)- \nu_{0}\left(f\right) - \nu_0'\left(f\right)\right| \leq \frac{K\left(n\right)}{N} \nu\left(f^2\right)^{\frac 12}.
\end{equation}
Substituting for $f$ the appropriate functions from our hypothesis, and using the independence of the last coordinate from the first $N-1$ in the $\nu_0$ measure, we see that the second term on the left vanishes in all cases.  It therefore suffices to compute 
\begin{equation}
\nu_{0}' \left(\prod_{r\in S} \left(\varepsilon^{r}- \varepsilon^{\zeta\left(r\right)}\right) \prod_{r'=1}^p \varepsilon^{\eta\left(r'\right)}f^{-}\right) \text{ and } \nu_0'\left(\left(\prod_{l=1}^p \varepsilon^{\eta\left(l\right)}-q_p\right)f^{-}\right).
\end{equation}

Applying Formula 2.169 from \cite{Talagrand-book},
\begin{multline}\label{Eq:cavitation}
\nu_{0}' \left(\prod_{r\in S} \left(\varepsilon^{r}- \varepsilon^{\zeta\left(r\right)}\right) \prod_{l=1}^p \varepsilon^{\eta\left(l\right)}f^{-}\right) =\\
 \beta^2 \sum_{\{\ell, \ell'\} \subseteq [n]}\nu_{0}\left(\prod_{r\in S} \left(\varepsilon^{r}- \varepsilon^{\zeta\left(r\right)}\right) \prod_{l=1}^p \varepsilon^{\eta\left(l\right)}\varepsilon^{\ell}\varepsilon^{\ell'}\right) \nu_{0}\left(f^-\left(R_{\ell, \ell'}^--q_2\right)\right)\\
- \beta^2 n \sum_{1 \leq \ell \leq [n]}\nu_{0}\left(\prod_{r\in S} \left(\varepsilon^{r}- \varepsilon^{\zeta\left(r\right)}\right) \prod_{1 \leq l \leq p} \varepsilon^{\eta\left(l\right)} \varepsilon^{\ell} \varepsilon^{n+1}\right) \nu_{0}\left(f^-\left(R_{\ell, n+1}^--q_2\right)\right)\\
+  \beta^2 \frac{n\left(n+1\right)}{2} \nu_{0}\left(\prod_{r\in S} \left(\varepsilon^{r}- \varepsilon^{\zeta\left(r\right)}\right) \prod_{l=1}^p \varepsilon^{\eta\left(l\right)} \varepsilon^{n+1} \varepsilon^{n+2}\right) \nu_{0}\left(f^-\left(R_{n+1, n+2}^--q_2\right)\right)
 \end{multline}
and similarly for $\prod_{l=1}^p \varepsilon^{\eta\left(l\right)}-q_p$.

Notice first that
\begin{equation}
\nu_{0}\left(\prod_{r\in S} \left(\varepsilon^{r}- \varepsilon^{\zeta\left(r\right)}\right) \prod_{l=1}^p \varepsilon^{\eta\left(l\right)}\varepsilon^{\ell}\varepsilon^{\ell'}\right) = 0
\end{equation}
whenever there is a factor $(\varepsilon^r- \varepsilon^{\zeta(r)})$ 'free' in the sense that $\{\ell, \ell'\}\cap  \{r,\zeta(r)\}= \varnothing$ or when $\{\ell, \ell'\} = \{r,\zeta(r)\}$.   When $|S| \geq 3$ this is always true, so the first assertion of the lemma follows.

Consider the case $|S|=2$.  For definiteness let us take $S=\{r_1, r_2\}$.
After the dust settles the only contributions come from the cases
\begin{equation}
\{\ell , \ell'\} \in \left\{\{r_1,\zeta(r_2)\},\{r_2, \zeta(r_1)\}\right\}.
\end{equation}
This leads to the expression
\begin{multline}
\nu_{0}' \left(\left(\varepsilon^{r_1}- \varepsilon^{\zeta\left(r_1\right)}\right)\left(\varepsilon^{r_2}- \varepsilon^{\zeta\left(r_2\right)}\right) \prod_{l=1}^p \varepsilon^{\eta\left(i\right)}f^{-}\right) =\\
\beta^2 \nu_{0}\left(\left(\varepsilon^{i_1}- \varepsilon^{\zeta\left(i_1\right)}\right)\left(\varepsilon^{r_2}- \varepsilon^{\zeta\left(r_2\right)}\right) \prod_{l=1}^p \varepsilon^{\eta\left(l\right)}\varepsilon^{r_1}\varepsilon^{r_2}\right) \times\\
\nu_{0}\left(f^-\left(R_{r_1, r_2}^--R_{r_1, \zeta\left(r_2\right)}^- - R_{\zeta\left(r_1\right), r_2}^-+R_{\zeta\left(r_1\right), \zeta\left(r_2\right)}^-\right)\right).
\end{multline}
If we calculate the first factor in terms of the values $\{q_s\}_{s \in \N}$ the assertion for $|S|=2$ follows as well.

Consider next the case $|S|=1$.  Letting $S=\{r_0\}$ for definiteness, we can, by replica symmetry, conclude that the summands in \eqref{Eq:cavitation} for which $\{\ell, \ell' \} \cap \{r_0, \zeta\left(r_0\right)\}= \varnothing$ when $\{\ell, \ell'\} = \{r, \zeta\left(r\right)\}$ all vanish.
Since
\begin{equation}
\nu_{0}\left(\left(\varepsilon^{r_0}- \varepsilon^{\zeta\left(r_0\right)}\right) \prod_{i=1}^p \varepsilon^{\eta\left(l\right)}\varepsilon^{r_0}\varepsilon^{\ell}\right) = - \nu_{0}\left(\left(\varepsilon^{r_0}- \varepsilon^{\zeta\left(r_0\right)}\right) \prod_{l=1}^p \varepsilon^{\eta\left(l\right)}\varepsilon^{\zeta(r_0)}\varepsilon^{\ell}\right),
\end{equation}
we have
\begin{multline}\label{Eq:cavitation2}
\nu_{0}' \left(\left(\varepsilon^{r_0}- \varepsilon^{\zeta\left(r_0\right)}\right) \prod_{l=1}^p \varepsilon^{\eta\left(l\right)}f^{-}\right) =\\
 \beta^2 \sum_{1 \leq \ell \leq p}\nu_{0}\left(\left(\varepsilon^{r_0}- \varepsilon^{\zeta\left(r_0\right)}\right) \prod_{l=1}^p \varepsilon^{\eta\left(l\right)}\varepsilon^{r_0}\varepsilon^{\eta\left(\ell\right)}\right) \nu_{0}\left(f^-\left(R_{r_0, \eta\left(\ell\right)}^--R_{\zeta\left(r_0\right), \eta\left(\ell\right)}^-\right)\right)\\
+ \beta^2 \sum_{ \ell \notin\{r_0,\zeta\left(r_0\right)\}\cup \eta\left([p]\right)}\nu_{0}\left(\left(\varepsilon^{r_0}- \varepsilon^{\zeta\left(r_0\right)}\right) \prod_{l=1}^p \varepsilon^{\eta\left(l\right)}\varepsilon^{r_0}\varepsilon^{\ell}\right) \nu_{0}\left(f^-\left(R_{r_0, \ell}^--R_{\zeta\left(r_0\right), \ell}^-\right)\right)\\
-  \beta^2 n \nu_{0}\left(\left(\varepsilon^{r_0}- \varepsilon^{\zeta\left(r_0\right)}\right) \prod_{l=1}^p \varepsilon^{\eta\left(l\right)} \varepsilon^{r_0} \varepsilon^{n+1}\right) \nu_{0}\left(f^-\left(R_{r_0, n+1}^--R_{\zeta\left(i_0\right), n+1}^-\right)\right).
 \end{multline}
 The assertion now follows via calculation of the factors involving $\varepsilon$'s in terms of the values $\{q_s\}_{s \in \N}$.
 
In the final case, no \textit{a priori} cancelations can be made.  Our assertion follows directly from a calculation of the resulting factors involving the $\varepsilon$'s.
\end{proof}

Our strategy is to use the previous lemma and an induction argument to calculate moments to within an error which is of higher order than the moment under consideration.  This entails getting bounds on the error term of the previous lemma, when the function $f$ is a monomial in the truncated multi-overlaps.  The next lemma follows (modulo a bit of work) from the exponential bounds for overlaps proved in 
Section 2.5 of \cite{Talagrand-book}.  We shall provide a proof in the Appendix in order to keep this paper self-contained.  
Let us introduce here a notation that is pervasive below.  Let 
\begin{equation}
T_{S,p}^-= \sum_{i=1}^{N-1} \prod_{r \in S}\left(\sigma_i^r- \langle\sigma_i\rangle\right) \langle\sigma_i\rangle^p
\end{equation}
\begin{lemma}
\label{L:Mom-Bound}
For each $n, m \in \mathbb N$, consider the joint collection $\{T_{S,p}\}_{S \subset [n],0 \leq p \leq m}$.  Then for any fixed collection of integers $\{k\left(S,p\right)\}_{S \subset [n], 0 \leq p \leq m}$
\begin{equation}
\left|\nu \left(\prod_{S \subset [n]} \prod_{0 \leq p \leq m} T_{S,p}^{k\left(S,p\right)}\right)\right| \leq \frac{C}{N^{\frac{k}{2}}}
\end{equation}
where $k= \sum_{S, p} k\left(S,p\right)$.
Similarly,
for any fixed collection of integers $\{k\left(S,p\right)\}_{S \subset [n], 0 \leq p \leq m}$
\begin{equation}
\left|\nu \left(\prod_{S \subset [n]} \prod_{0 \leq p \leq m} \left(T_{S,p}^-\right)^{k\left(S,p\right)}\right)\right| \leq \frac{C}{N^{\frac{k}{2}}}
\end{equation}
where $k= \sum_{S, p} k\left(S,p\right)$.
\end{lemma}

We proceed to identify the covariance structure of the process $\{T_{S,p}\}$ in a series of lemmas.  We begin with observations which are the basis of the calculations below.  For ease of notation, let us restrict to the case $|S|,|\tilde S| \geq 1$.  The remaining cases are similar.

Using the symmetry of $\nu$,
\begin{equation}
\nu\left(T_{S,p}T_{\tilde{S}, \tilde{p}}\right)  = \nu\left(\prod_{r \in S}\left(\varepsilon^r- \varepsilon^{\zeta\left(r\right)}\right)\prod_{1 \leq l \leq  p} \varepsilon^{\eta\left(l\right)} T_{\tilde{S}, \tilde{p}}\right)
\end{equation}
where $\zeta,\eta$ are injective integer valued functions as defined above Lemma \ref{L:cavitation}.

Note that
\begin{equation}
T_{\tilde{S}, \tilde{p}} = T_{\tilde{S}, \tilde{p}}^{-} + \frac{1}{N} \prod_{i \in \tilde{S}}\left(\varepsilon^i- \varepsilon^{\alpha\left(i\right)}\right)\prod_{1 \leq l \leq \tilde p} \varepsilon^{\gamma\left(l\right)},
\end{equation}
where $\alpha: \tilde{S} \rightarrow \N$ and $\gamma: [\tilde p] \rightarrow \N$ are injective with ranges disjoint from $\{S, \tilde S, \zeta(S), \eta([p])\}$ and each other.
We have
\begin{multline}
\label{Eq:nu_0-1}
\nu\left(\prod_{r \in S}\left(\varepsilon^r- \varepsilon^{\zeta\left(r\right)}\right)\prod_{1 \leq l \leq p} \varepsilon^{\eta\left(l\right)}T_{\tilde{S}, \tilde{p}}\right)  = \\
\frac{1}{N} \nu\left(\prod_{r \in S}\left(\varepsilon^r- \varepsilon^{\zeta\left(r\right)}\right)\prod_{1 \leq l \leq p} \varepsilon^{\eta\left(l\right)}  \prod_{r \in \tilde{S}}\left(\varepsilon^r- \varepsilon^{\alpha\left(r\right)}\right)\prod_{1 \leq l\leq \tilde{p}} \varepsilon^{\gamma\left(l\right)}\right) \\+ \nu\left(\prod_{r \in S}\left(\varepsilon^r- \varepsilon^{\zeta\left(r\right)}\right)\prod_{1 \leq l \leq p} \epsilon^{\eta\left(l\right)} T_{\tilde{S}, \tilde{p}}^-\right)
\end{multline}

\begin{lemma}
\label{L:Ind}
Suppose $n, m \in \N$ is fixed.  Let $S, \tilde S \subseteq [n]$ so that $S \neq \tilde S$.  Let $0 \leq p, \tilde p \leq  m$.  Then
\begin{equation}
\nu\left(T_{S,p}T_{\tilde S, \tilde p}\right) = O\left(3\right)
\end{equation} 

More generally,
let $\mathcal G$ be a monomial in the truncated multi-overlaps $\{T_{S',p'}\}_{S' \subseteq [n], p' \in [m]}$ with total degree $k$.  If $S_0 \neq S'$ for all factors of $\mathcal G$, then
\begin{equation}
\nu\left(T_{S_0,p_0}\mathcal G\right) = O\left(k+2\right)
\end{equation}
for any $p_0 \in \N$.
\end{lemma}
 
\begin{proof}
We shall only prove the statement for covariances.  The more general statement follows from similar, but more involved computations.  For notational convenience, we restrict attention to the case that $|S|, |\tilde S| \geq 1$, the remaining cases following similar arguments.
Consider first the term
\begin{equation}
 \frac{1}{N} \nu\left(\prod_{r \in S}\left(\varepsilon^r- \varepsilon^{\zeta\left(r\right)}\right)\prod_{1 \leq l \leq p} \varepsilon^{\eta\left(l\right)}  \prod_{r \in \tilde{S}}\left(\varepsilon^r- \varepsilon^{\alpha\left(r\right)}\right)\prod_{1 \leq l \leq \tilde{p}} \varepsilon^{\gamma\left(l\right)}\right)
 \end{equation}
from \eqref{Eq:nu_0-1}.  Using Lemma \ref{L:Cav-Bound} to approximate this term via the cavity method, notice that the left endpoint of the interpolation vanishes whenever $S \varDelta \tilde S \neq \varnothing$ via symmetry.  Thus in all cases this term is $O(3)$.   

We claim that
\begin{equation}
\nu\left(\prod_{r \in S}\left(\varepsilon^r- \varepsilon^{\zeta\left(r\right)}\right)\prod_{1 \leq l \leq p} \varepsilon^{\eta\left(l\right)} T_{\tilde{S}, \tilde{p}}^-\right)  = O(3)
\end{equation}
as well.  Without loss of generality, we may assume $|S| \geq |\tilde{S}|$.  The first observation to make is that if $|S| \geq 3$ then Lemmas \ref{L:cavitation} and \ref{L:Mom-Bound} give the conclusion immediately.

The remaining cases follow from more delicate cancellations.  We may assume (with our initial assumption that neither $S$ nor $\tilde{S}$ is empty still in place) that $|S|=2$ and will denote $S=\{r_1,r_2\}$.  By Lemmas  \ref{L:Cav-Bound} and \ref{L:cavitation}
\begin{multline}
 \nu\left(\prod_{r \in S}\left(\varepsilon^r- \varepsilon^{\zeta\left(r\right)}\right)\prod_{1 \leq l \leq p} \varepsilon^{\eta\left(l\right)} T_{\tilde{S}, \tilde{p}}^-\right) = \\
\beta^2\varrho_p
\nu_{0}\left(T^-_{\tilde{S}, \tilde{p}}\left(R_{r_1, r_2}^--R_{r_1, \zeta\left(r_2\right)}^- - R_{\zeta\left(r_1\right), r_2}^-+R_{\zeta\left(r_1\right), \zeta\left(r_2\right)}^-\right)\right)
 + O\left(3\right)
\end{multline}
where we have again employed Lemma \ref{L:Mom-Bound} to bound the error term.
Applying Lemma \ref{L:Cav-Bound},
\begin{multline}
\nu_{0}\left(T^-_{\tilde{S}, \tilde{p}}\left(R_{r_1, r_2}^--R_{r_1, \zeta\left(r_2\right)}^- - R_{\zeta\left(r_1\right), r_2}^-+R_{\zeta\left(r_1\right), \zeta\left(r_2\right)}^-\right)\right) =\\
\nu\left(T_{\tilde{S}, \tilde{p}}\left(R_{r_1, r_2}-R_{r_1, \zeta\left(r_2\right)} - R_{\zeta\left(r_1\right), r_2}+R_{\zeta\left(r_1\right), \zeta\left(r_2\right)}\right)\right) + O(3)
\end{multline}
By assumption, $\tilde{S}$ can involve at most one of $\{r_1, r_2\}$,
it is easy to check that 
\begin{equation}
\nu\left(T_{\tilde{S}, \tilde{p}}\left(R_{r_1, r_2}-R_{r_1, \zeta\left(r_2\right)} - R_{\zeta\left(r_1\right), r_2}+R_{\zeta\left(r_1\right), \zeta\left(r_2\right)}\right)\right) = 0
\end{equation}
The remaining cases follow from similar considerations.
\end{proof}

The previous lemma will be key in diminishing the pain of determining non-trivial covariances.  The other helpful step was taken already by Talagrand, who computed the variances of $T_{\{1,2\}}, T_{1,1}$ and $T_{\varnothing, 2}$.  As should not be a surprise, these prove to be the most fundamental truncated overlaps.  To prove the remainder of Theorem \ref{T:Cov}, we must specialize the analysis to the various cases.  

\section{The Covariance Structure}
\label{four}
\noindent We begin this section with a small computation.  By employing the Lemma \ref{L:Cav-Bound} and computing the left endpoint expectation directly, we have the preliminary observation
\begin{multline}
\label{Eq:0-order}
\frac{1}{N} \nu\left(\prod_{r \in S}\left(\varepsilon^r- \varepsilon^{\zeta\left(r\right)}\right)\prod_{1 \leq l \leq p} \varepsilon^{\eta\left(l\right)}  \prod_{r \in S}\left(\varepsilon^r- \varepsilon^{\alpha\left(r\right)}\right)\prod_{1 \leq l \leq \tilde{p}} \varepsilon^{\gamma\left(l\right)}\right)= \\
\frac{1}{N}\sum_{\ell = 0}^{|S|} {|S| \choose \ell} \left(-1\right)^{\ell} q_{p+\tilde{p} + 2 \ell} + O(3)
\end{multline}

We begin with the simplest case:

\begin{lemma}
\label{L:|S|=3}
Suppose $n \in \N$ is fixed.
Let $S\subseteq [n]$ and $p, \tilde p \in \N$ be fixed.  If $|S| \geq 3$, then
\begin{equation}
\nu\left(T_{S,p}T_{S, \tilde{p}}\right) = \frac{1}{N}\sum_{\ell = 0}^{|S|} {|S| \choose \ell} \left(-1\right)^{\ell} q_{p+\tilde{p} + 2 \ell} + O(3)
\end{equation}
\end{lemma}

\begin{proof}
We shall use \eqref{Eq:nu_0-1}.  Since $|S| \geq 3$, Lemmas \ref{L:cavitation} and \ref{L:Mom-Bound} imply
\begin{equation}
\nu\left(\prod_{i \in S}\left(\varepsilon^r- \varepsilon^{\zeta\left(r\right)}\right)\prod_{1 \leq l \leq  p} \varepsilon^{\eta\left(l\right)} T_{\tilde{S}, \tilde{p}}^-\right)  = O(3)
\end{equation}
The assertion then follows immediately from \eqref{Eq:0-order}.
\end{proof}

\begin{lemma}
\label{L:|S|=2}
Suppose $n \in \N$ is fixed.
Let $S\subseteq [n]$ and $p, \tilde p \in \N$ be fixed.  If $|S| =2$, then
\begin{align*}
\nu\left(T_{S,p}T_{S, \tilde{p}}\right) = &\frac{1}{N} \varrho_{p+\tilde p} + \frac{\beta^2}{N}\varrho_{p} \varrho_{\tilde p} + \frac{\beta^4}{N}\varrho_{p}  \varrho_{\tilde p}  \frac{\varrho_{0}}{1- \beta^2\varrho_{0}}.
 \end{align*}
\end{lemma}

\begin{proof}
From \eqref{Eq:nu_0-1} and \eqref{Eq:0-order}, it is enough to identify the value of 
\begin{equation}
\nu\left(\prod_{r\in S} \left(\varepsilon^{r}- \varepsilon^{\zeta\left(r\right)}\right) \prod_{1 \leq l \leq p} \varepsilon^{\eta\left(l\right)}T_{S, \tilde p}^{-}\right).
\end{equation}
Let $S=\{r_1, r_2\}$ for definiteness.  Note that 
\begin{equation}
R_{r_1, r_2}-R_{r_1, \zeta\left(r_2\right)} - R_{\zeta\left(r_1\right), r_2}+R_{\zeta\left(r_1\right), \zeta\left(r_2\right)}=T_{r_1, r_2}-T_{r_1, \zeta\left(r_2\right)} - T_{\zeta\left(r_1\right), r_2}+T_{\zeta\left(r_1\right), \zeta\left(r_2\right)}
\end{equation}
With this identity we have, by applying Lemmas \ref{L:Cav-Bound}, \ref{L:cavitation}, \ref{L:Mom-Bound} and then Lemma \ref{L:Cav-Bound} once again
\begin{multline}
\nu\left(\prod_{r\in S} \left(\varepsilon^{r}- \varepsilon^{\zeta\left(r\right)}\right) \prod_{1 \leq l \leq p} \varepsilon^{\eta \left(l\right)}T^{-}_{ S, \tilde{p}} \right)= \\
\beta^2\varrho_p
\nu_{0}\left(T^-_{S, \tilde{p}} \left(R_{r_1, r_2}^--R_{r_1, \zeta\left(r_2\right)}^- - R_{\zeta\left(r_1\right), r_2}^-+R_{\zeta\left(r_1\right), \zeta\left(r_2\right)}^-\right)\right) + O\left(3\right)\\
=
\beta^2\varrho_{p}\nu\left(T_{S,\tilde{p}} \left(T_{r_1, r_2}-T_{r_1, \zeta\left(r_2\right)} - T_{\zeta\left(r_1\right), r_2}+T_{\zeta\left(r_1\right), \zeta\left(r_2\right)}\right)\right) + O\left(3\right).
\end{multline}
Now Lemma \ref{L:Ind} implies
\begin{equation}
\nu\left(T_{S, \tilde{p}} \left(T_{r_1, r_2}-T_{r_1, \zeta\left(r_2\right)} - T_{\zeta\left(r_1\right), r_2}+T_{\zeta\left(r_1\right), \zeta\left(r_2\right)}\right)\right)= \nu\left(T_{S,\tilde{p}} T_{r_1, r_2}\right) + O(3).
\end{equation}

To summarize, we have shown that
\begin{equation}
\nu\left(T_{S,p}T_{S, \tilde{p}}\right) =\frac{1}{N} \varrho_{p+\tilde p} + \beta^2\varrho_{p}\nu\left(T_{S,\tilde p}T_S\right) + O(3).
\end{equation}
Since $p$ and $\tilde p$ were chosen arbitrarily in the preceding, specializing to $\tilde p=0$ we may employ this approximate identity to 
see that 
\begin{equation}
\nu\left(T_{S, \tilde p}T_S\right)=\frac{1}{N} \varrho_{\tilde p} +\beta^2\varrho_{\tilde p}\nu\left(T_S^2\right) + O(3).
\end{equation}

To finish, we rely on Theorem \ref{T:CLT-Tal} which tells us that
\begin{equation}
\nu\left(T_S^2\right) =  \frac{\varrho_{0}}{1- \beta^2\varrho_{0}} +O(3).
\end{equation}
Combing these approximations together gives the result.
\end{proof}
For the next two Lemmas (the cases $|S| \leq 1$), it is convenient (not to mention more compact and instructive) to give approximate identities for the covariances rather than full formulas.  One can easily construct the corresponding formulas iteratively.

\begin{lemma}
\label{L:|S|=1}.
Let $r_1, p, \tilde p \in \N$ be fixed.  Then
\begin{multline}
 \label{Eq:S=1Id}
\nu\left(T_{r_1, p}T_{r_1,\tilde{p}}\right)=  \frac{1}{N}\pi_{p+\tilde{p}-1}+
\beta^2\left[p \pi_{p-2} - (p+2)\pi_p\right] \nu\left(T_{r_1, \tilde p} T_{r_1, 1}\right)\\
 + \beta^2\tilde p\pi_{p} \nu\left(T_{\{1,2\}, \tilde{p}-1}T_{\{1,2\}}\right)
+ O(3)
\end{multline}
where we interpret a term with a coefficient $0$ as evaluating to $0$.
\end{lemma}

\begin{proof}
To fix notation, let $S=\{r_1\}$.  Let $\zeta,\eta, \alpha$ and $\gamma$ be integer valued functions as in \eqref{Eq:nu_0-1}.
By Lemma \ref{L:Cav-Bound}
\begin{equation}
\frac{1}{N} \nu\left(\left(\varepsilon^{r_1}- \varepsilon^{\zeta\left(r_1\right)}\right)\prod_{1 \leq l \leq p} \varepsilon^{\eta \left(l\right)}  \left(\varepsilon^{r_1}- \varepsilon^{\alpha\left(r_1\right)}\right)\prod_{1 \leq l \leq \tilde{p}} \varepsilon^{\gamma\left(l\right)}\right)= \frac{\pi_{p+\tilde p-1}}{N}
\end{equation}
Thus, from \eqref{Eq:nu_0-1}, it is enough to identify the approximate value of 
\begin{equation}
\nu\left(\left(\varepsilon^{r_1}- \varepsilon^{\zeta\left(r_1\right)}\right) \prod_{1\leq l \leq p} \varepsilon^{\eta\left(l\right)}T_{r_1,\tilde{p}}^{-}\right).
\end{equation}

As usual, we apply Lemmas \ref{L:Cav-Bound}, \ref{L:cavitation} and \ref{L:Mom-Bound}.  Note that the Lemmas \ref{L:Cav-Bound} and \ref{L:Mom-Bound} imply that for any triple $a, b, c \in \N$,
\begin{equation}
 \nu_{0}\left(T_{r_1, \tilde{p}}^-\left(R_{a, b}^--R_{c, b}^-\right)\right) =   \nu\left(T_{r_1, \tilde{p}}\left(R_{a, b}-R_{c, b}\right)\right) + O(3).
 \end{equation}
Further, if $T_{r_1, \tilde{p}}$ does not depend on the replicas $b, c$
\begin{equation}
\nu\left(T_{r_1, \tilde{p}}\left(R_{a, b}-R_{c, b}\right)\right)= \nu\left(T_{r_1, \tilde{p}}T_{r_1,1}\right).
\end{equation}
The resulting symmetry allows us collect terms to get
\begin{multline}
\nu\left(\left(\varepsilon^{r_1}- \varepsilon^{\zeta\left(r_1\right)}\right) \prod_{1 \leq l \leq p} \varepsilon^{\eta\left(l\right)}T_{r_1,\tilde{p}}\right)= \beta^2\left[p\pi_{p-2} -(p+\tilde p + 3)\pi_{p}\right] \nu\left(T_{r_1, \tilde p} T_{r_1,1}\right) \\
+ \beta^2 p \pi_{p}  \nu\left(T_{r_1, \tilde p}(\alpha(1)) \left(R_{r_1,\alpha(1)}- R_{\zeta(r_1), \alpha(1)}\right)\right)\\
+ \beta^2\tilde p\pi_{p} \nu\left(T_{r_1, \tilde p}(\gamma(1)) \left(R_{r_1,\gamma(1)}- R_{\zeta(r_1), \gamma(1)}\right)\right) + O(3).
\end{multline}

To compute the last two summand we introduce the identities 
\begin{align*}
R_{a,b}=& T_{\{a,b\}} + T_{a,1}+ T_{b,1} + T + q_2\\
R_{r_1,\gamma(1)}- R_{\zeta(r_1), \gamma(1)} =& T_{\{r_1, \gamma(1)\}}- T_{\{\zeta(r_1), \gamma(1)\}}+ T_{r_1,1} - T_{\zeta(r_1),1}\\
R_{r_1,\alpha(r_1)}- R_{\zeta(r_1), \alpha(r_1)} =&T_{\{r_1, \alpha(r_1)\}} - T_{\{\zeta(r_1), \alpha(r_1)\}}+T_{r_1, 1}-T_{\zeta(r_1), 1}
\end{align*}
Further, we have
\begin{align*}
T_{r_1, \tilde p}(\alpha(1))=& T_{r_1, \tilde p}-T_{\alpha(1), \tilde p}\\
T_{r_1,\tilde p}(\gamma(1))=&T_{r_1, \tilde p}(\gamma(1))- T_{r_1, \tilde p}(z_0) + T_{r_1, \tilde p}(z_0)= T_{\{r_1, \gamma(1)\},\tilde p-1} + T_{r_1,\tilde{p}}
\end{align*}
where the index $z_0$ is disjoint from those already under consideration and the second identity only makes sense if $\tilde p >0$.

Let us now collect terms, applying Lemma \ref{L:Ind} whenever possible.  We have
\begin{align*}
\nu\left(T_{r_1, \tilde p}(\alpha(1)) \left(R_{r_1,\alpha(1)}- R_{\zeta(r_1), \alpha(1)}\right)\right)=&\nu\left(T_{r_1, \tilde p} T_{r_1, 1}\right) + O(3)\\
\nu\left(T_{r_1, \tilde p}(\gamma(1)) \left(R_{r_1,\gamma(1)}- R_{\zeta(r_1), \gamma(1)}\right)\right)=&  \nu\left(T_{\{1, 2\}, \tilde p-1} T_{\{1,2\}}\right) + \nu\left(T_{r_1,\tilde p}T_{\{r_1\}, 1}\right) + O(3)
\end{align*}
This gives us the approximate identity
\begin{multline}
\nu\left(T_{r_1, p}T_{r_1,\tilde{p}}\right)=  \frac{1}{N}\pi_{p+\tilde{p}-1}+
\beta^2\left[p \pi_{p-2} - 3\pi_p\right] \nu\left(T_{r_1, \tilde p} T_{r_1, 1}\right)\\
 + \beta^2\tilde p\pi_{p} \nu\left(T_{\{1,2\}, \tilde{p}-1}T_{\{1,2\}}\right)
+ O(3)
\end{multline}
where we interpret a term with a coefficient $0$ as evaluating to $0$.
\end{proof}

\begin{lemma}
\label{L:|S|=0}
Suppose $n \in \N$ is fixed.
Let $p, \tilde p \in \N$ be fixed.  We have
\begin{multline}
\nu\left(T_{\varnothing, p}T_{\varnothing, \tilde p}\right)= \frac{1}{N}\left(q_{p+\tilde p}-q_pq_{\tilde p}\right) +  \beta^2 \lambda _{p} \nu\left(T_{\varnothing, \tilde p}T_{\varnothing, 2}\right) 
+\beta^2\kappa_{p, \tilde p} \nu\left(T_{1, \tilde p-1}T_{1,1}\right) \\
+\beta^2 \omega_{p, \tilde p}\nu\left(T_{\{1,2\}, \tilde p-2}T_{\{1,2\}}\right) + O(3).
\end{multline}
\end{lemma}

\begin{proof}
Let $\eta$ and $\gamma$ denote functions as in \eqref{Eq:nu_0-1} which map  $[p]$ and $[\tilde p]$, respectively, into $\N$.
By Lemma \ref{L:Cav-Bound}
\begin{equation}
\frac{1}{N}\nu\left(\left(\prod_{l=1}^p \varepsilon^{\eta\left(l\right)} - q_p\right) \left(\prod_{l=1}^{\tilde p}\varepsilon^{\gamma \left(l\right)} - q_{\tilde p}\right)\right) = \frac{1}{N}\left(q_{p+\tilde p}-q_pq_{\tilde p}\right) + O(3).
\end{equation}
From \eqref{Eq:nu_0-1} and \eqref{Eq:0-order}, it is enough to identify the value of 
\begin{equation}
\nu\left(\left(\prod_{l=1}^p \varepsilon^{\eta\left(l\right)} - q_p\right) T_{\varnothing,\tilde{p}}^{-}\right).
\end{equation}

In the following calculation, we write 
\begin{equation}
T_{\varnothing, \tilde p}= \frac{1}{N} \sum_{i=1}^N  \sigma_i^{\gamma(1)} \cdots \sigma_i^{\gamma(\tilde p)}.
\end{equation}

Lemmas \ref{L:Cav-Bound} and \ref{L:Mom-Bound} imply that for any pair $(a, b)$,
\begin{equation}
\label{Eq:return}
 \nu_{0}\left(T_{\varnothing, \tilde{p}}^- \left(R_{a, b}^--q_2\right)\right) =   \nu\left(T_{\varnothing, \tilde{p}}\left(R_{a, b}-q_2\right)\right) + O(3).
 \end{equation}
Recall Lemma \ref{L:cavitation} in the present case.  The symmetry resulting from \eqref{Eq:return} allows us to write
\begin{multline}
\nu\left(\left(\prod_{l=1}^p \varepsilon^{\eta\left(l\right)}- q_p\right)T_{\varnothing, \tilde p}^{-}\right)=
\beta^2 \left(q_{p-2}- q_pq_2\right)\frac{p(p-1)}{2}\nu\left(T_{\varnothing, \tilde p} \left(R_{\eta(1), \eta(2)}-q_2\right)\right)\\ 
+\beta^2 \left(q_p- q_p q_2\right)p \tilde p\nu\left(T_{\varnothing, \tilde p}(\gamma(1))\left(R_{\eta(1), \gamma(1)}-q_2\right)\right)\\
+\beta^2 \left(q_{p+2}-q_pq_2\right) \frac{\tilde{p} (\tilde{p}-1)}{2}\nu\left(T_{\varnothing, \tilde p}(\gamma(1), \gamma(2))\left(R_{\gamma(1) \gamma(2)}-q_2\right)\right)\\
- \beta^2 \left(q_p-q_pq_2\right)(p+\tilde p)p  \nu\left(T_{\varnothing, \tilde p} \left(R_{\eta(1), z}-q_2\right)\right)\\
- \beta^2 \left(q_{p+2}-q_pq_2\right)(p+\tilde p)\tilde p \nu\left(T_{\varnothing, \tilde p}(\gamma(1))\left(R_{\gamma(1), z}-q_2\right)\right)\\
+  \beta^2 \frac{(p+\tilde p)\left(p+\tilde p+1\right)}{2} \left(q_{p+2}-q_p q_2\right)\nu\left(T_{\varnothing, \tilde p}\left(R_{z, z+1}-q_2\right)\right).
\end{multline}

For each pair $\{a,b\}$, we may write
\begin{equation}
R_{a,b}-q_2 = T_{\{a,b\}} + T_{a, 1}+ T_{b,1} + T_{\varnothing, 2}.
\end{equation}
For the first, fourth and sixth terms, Lemma \ref{L:Ind} and this expansion allow us to replace $\nu\left(T_{\varnothing, \tilde p} \left(R_{a, b}-q_2\right)\right)$ by $\nu\left(T_{\varnothing, \tilde p} T_{\varnothing, 2}\right) + O(3)$

To handle the second and fifth terms we introduce a new replica $z$ which replaces $\gamma(1)$ and write
\begin{align*}
T_{\varnothing, \tilde p}(\gamma(1))= &T_{\varnothing, \tilde p}(\gamma(1))- T_{\varnothing, \tilde p}(z) +T_{\varnothing, \tilde p}(z)\\
=& T_{\gamma(1), \tilde p-1} + T_{\varnothing, \tilde p}
\end{align*}
Applying Lemma \ref{L:Ind} once more gives, via symmetry, 
\begin{equation}
\nu\left(T_{\varnothing, \tilde p}(\gamma(1)) \left(R_{\gamma(1), b}-q_2\right)\right)=\nu\left(T_{1, \tilde p-1} T_{1, 1}\right) + \nu\left(T_{\varnothing, \tilde p} T_{\varnothing, 2}\right) + O(3).
\end{equation}

To handle the third term we make a similar expansion of $T_{\varnothing, \tilde p}(\gamma(1), \gamma(2))$.  We find
\begin{multline}
\nu\left(T_{\varnothing, \tilde p}(\gamma(1), \gamma(2))\left(R_{\gamma(1) \gamma(2)}-q_2\right)\right)=\\ 
\nu\left(T_{\{1,2\}, \tilde p-2}T_{\{1, 2\}}\right)+2 \nu\left(T_{1, \tilde p-1}T_{1, 1}\right)+ \nu\left(T_{\varnothing, \tilde p}T_{\varnothing , 2}\right) + O(3)
\end{multline}
All in all, we have
\begin{multline}
\nu\left(T_{\varnothing, p}T_{\varnothing, \tilde p}\right)=  \frac{1}{N}\left(q_{p+\tilde p}-q_pq_{\tilde p}\right) + \beta^2 \lambda _{p} \nu\left(T_{\varnothing, \tilde p}T_{\varnothing, 2}\right) 
+\beta^2\kappa_{p, \tilde p} \nu\left(T_{1, \tilde p-1}T_{1,1}\right) \\
+\beta^2 \omega_{p, \tilde p}\nu\left(T_{\{1,2\}, \tilde p-2}T_{\{1,2\}}\right) + O(3).
\end{multline}
\end{proof}

\section{The Identification of Moments}
Finally we are in a position to compute the joint moments of the truncated overlaps $\{T_{S, p}\}$.  As usual, we work from the easiest to hardest case.

\begin{theorem}[Joint Moments:$|S| \geq 3$]
\label{T:|S|geq3}
Let $m, n \in \N$ be fixed.  Consider the collection of subsets $S \subseteq [n]$ and the associated family $\{T_{S,p}\}_{S \subseteq [n], 0 \leq p \leq m}$.  Let $\{k\left(S,p\right)\}_{S \subseteq [n], p \leq m}$ be a collection of non-negative integers such that $k(S, p) = 0$ if $|S| \leq 2$.  Let $k= \sum_{p=1}^m k\left(S,p\right)$.  Then
\begin{equation}
\left|\nu\left(\prod_{\stackrel{S \subseteq [n]}{0 \leq p \leq m}}T_{S,p}^{k\left(S,p\right)}\right) - N^{- \ffrac{k}{2}} \prod_{S\subset [n]}\E \left[\prod_{p=0}^m Y_{S,p}^{k\left(S,p\right)}\right]\right| \leq \frac{C}{N^{\frac{k+1}{2}}}.
\end{equation}
\end{theorem}
\begin{proof}
We shall prove this statement by induction on $k$, the case $k=1$ is an easy computation based on our current knowledge (Lemma \ref{L:Cav-Bound}), and $k=2$ was proved as part of Theorem \ref{T:Cov} in Section \ref{four}.  Let us distinguish one particular pair $(S_0, p_0)$ such that $k(S_0, p_0) \neq 0$.  For notational convenience let
\begin{equation}
\mathcal{F}= T_{S_0,p_0}^{k(S_0, p_0)-1} \prod_{(S,p) \neq (S_0,p_0)} T_{S,p}^{k(S,p)}.
\end{equation}
By symmetry, we have
\begin{equation}
\nu\left(\prod_{\stackrel{S \subseteq [n]}{0 \leq p \leq m}} T_{S,p}^{\eta \left(S,p\right)}\right) = \nu\left(\prod_{r \in S_0} (\varepsilon^r- \varepsilon^{\zeta(r)}) \prod_{l=1}^{p_0} \varepsilon^{\eta(l)} \mathcal F\right).
\end{equation}

We expand each of the factors $T_{S', p'}$ in $\mathcal F$ around $T_{S',p'}^-$, with $\mathcal F^-$ denoting the constant order term.  Applying Lemma \ref{L:Mom-Bound} to estimate the resulting summands, we have
\begin{multline}
\nu\left(\prod_{\stackrel{S \subseteq [n]}{0 \leq p \leq m}} T_{S,p}^{k\left(S,p\right)}\right) = \nu\left(\prod_{r \in S_0} (\varepsilon^r- \varepsilon^{\zeta(r)}) \prod_{1 \leq l \leq p_0} \varepsilon^{\eta(l)} \mathcal F^-\right)\\
 + \left(k(S_0,p_0)-1\right)  \frac{1}{N} \nu\left(\prod_{r\in S_0} (\varepsilon^r- \varepsilon^{\zeta(r)}) \prod_{1 \leq l \leq p_0} \varepsilon^{\eta(l)} \prod_{r \in S_0} (\varepsilon^r- \varepsilon^{\alpha(r)}) \prod_{1 \leq l \leq p_0}\varepsilon^{\gamma(l)}\frac{\mathcal F^-}{T_{S,p}^-}\right)\\
+ \sum_{(S,p) \neq (S_0, p_0)} k(S,p)  \frac{1}{N} \nu\left(\prod_{r\in S_0} (\varepsilon^r- \varepsilon^{\zeta(r)}) \prod_{1 \leq l \leq p_0} \varepsilon^{\eta(l)} \prod_{r \in S_0} (\varepsilon^r- \varepsilon^{\alpha(r)}) \prod_{1 \leq l \leq p_0}\varepsilon^{\gamma(l)}\frac{\mathcal F^-}{T_{S,p}^-}\right)\\
+ O(k+1).
\end{multline} 

Now Lemma \ref{L:Mom-Bound} implies
\begin{align*}
&\nu\left(\left[\mathcal F^-\right]^2\right)^{\frac 12} = O(k-1)\\
&\nu\left(\left[\frac{\mathcal F^-}{T_{S,p}^-}\right]^2\right)^{\frac 12} = O(k-2)
\end{align*}
when $k(S,p) \neq 0$.  Therefore, via Lemmas \ref{L:Cav-Bound} and \ref{L:cavitation} we have
\begin{multline}
\nu\left(\prod_{\stackrel{S \subseteq [n]}{0 \leq p \leq m}} T_{S,p}^{k\left(S,p\right)}\right) =\frac{1}{N}\left(k(S_0,p_0)-1\right)A_{|S_0|}(p_0,p_0) \nu \left(\frac{\mathcal F}{T_{S_0,p_0}}\right)\\+ \frac{1}{N}\sum_{(S,p) \neq (S_0,p_0)} k(S,p)\delta_{S, S_0}A_{|S_0|}(p,p_0) \frac{1}{N}\nu \left(\frac{\mathcal F}{T_{S,p}}\right)\\
+ O(k+1)
\end{multline} 
where we interpret a term in which a coefficient evaluates to $0$ as being $0$.
Expressed another way
\begin{multline}
\nu\left(\prod_{\stackrel{S \subseteq [n]}{0 \leq p \leq m}} T_{S,p}^{k\left(S,p\right)}\right) =\frac{1}{N}\sum_{\stackrel{S \subseteq [n]}{0 \leq p \leq m}}\E\left[Y_{S,p} Y_{S_0,p_0}\right] \nu \left(\frac{\partial}{\partial T_{S,p}}\mathcal F\right) + O(k+1).
\end{multline} 
From the induction hypothesis,
\begin{equation}
 \nu \left(\frac{\partial}{\partial T_{S,p}}\mathcal F\right) = \E\left[\frac{\partial}{\partial Y_{S,p}}\left(Y_{S_0,p_0}^{k(S_0,p_0)-1}\prod_{(S,p) \neq (S_0,p_0)}Y_{S,p}^{k(S,p)}\right)\right] + O(k-1).
\end{equation}
The theorem now follows easily from the assumption that the collection $\{Y_{S,p}\}$ has the specified Gaussian covariance structure.
\end{proof}

\begin{theorem}[Joint Moments:$|S| \geq 2$]
\label{T:|S|geq2}
Let $m, n \in \N$ be fixed.  Consider family $\{T_{S,p}\}_{S \subset [n], 0 \leq p \leq m}$ and let $\{k\left(S,p\right)\}_{S \subseteq [n],0 \leq  p \leq m}$ be a collection of non-negative integers such that $k(S, p) = 0$ if $|S| \leq 1$.  Let $k= \sum_{p=1}^m k\left(S,p\right)$.  Then
\begin{equation}
\left|\nu\left(\prod_{\stackrel{S \subseteq [n]}{0 \leq p \leq m} }T_{S,p}^{k\left(S,p\right)}\right) - N^{- \ffrac{k}{2}} \prod_{S \subseteq [n]}\E \left[\prod_{0 \leq p \leq m} Y_{S,p}^{k\left(S,p\right)}\right]\right| \leq \frac{C}{N^{\frac{k+1}{2}}}.
\end{equation}
\end{theorem}
\begin{proof}
Let
\begin{equation}
k_1=\sum_{\stackrel{S: |S| \geq 3}{0 \leq p \leq m}} k(S,p), \hspace{12pt} k_2= \sum_{\stackrel{|S|=2}{0 \leq p \leq m}} k(S,p))
\end{equation}
Let $k_2$ be fixed for the moment.  Then assuming that we have proved the statement when $k_1=0$, the case $k_1>0$ (and $k_2$ fixed) follows by an induction argument nearly exactly as in the previous theorem. Thus we may as well assume $k_1=0$.  

Let 
\begin{equation}
\tilde k_2 = \sum_{\stackrel{|S|=2}{p \geq 1}} k(S, p).
\end{equation}
We shall prove the statement (with $k_1 \equiv 0$) by induction on $\tilde k_2$.  The case $\tilde k_2 =0$ is part of Talgrand's CLT (Theorem \ref{T:CLT-Tal} in this paper).
The induction hypothesis for $\tilde k_2$ takes the following strong form.  Let us fix a collection of exponents $\{k(S,p)\}_{|S| = 2, 0 \leq p\leq m}$.  Suppose that $(S_0, p_0)$ is chosen so that $p_0$ is maximal with $k(S_0, p_0) \neq 0$.  
If there are multiple choices of $S_0$ we choose arbitrarily among those with $k(S_0,p_0)$ maximal as well.  As an induction hypothesis for $\tilde k_2$, let us suppose that we already have the statement for all choices of exponents $\{k^*(S, p)\}_{|S|= 2, 0\leq p \leq m}$ with $\{k^*(S,0)\}_{|S| = 2}$  arbitrary and  so that
\begin{equation}
k^*(S,p) \leq k(S,p),  \hspace{12pt} k^*(S_0,p_0) < k(S_0,p_0)
\end{equation}
for $|S| = 2$ and $p \geq 1$.
Thus Theorem \ref{T:CLT-Tal} amounts to the base case of our induction.

For notational convenience let
\begin{equation}
\mathcal{F}= T_{S_0,p_0}^{k(S_0, p_0)-1} \prod_{(S,p) \neq (S_0,p_0)} T_{S,p}^{k(S,p)}.
\end{equation}
By symmetry, we have
\begin{equation}
\nu\left(\prod_{\stackrel{S \subseteq [n]}{0 \leq p \leq m}} T_{S,p}^{k\left(S,p\right)}\right) = \nu\left(\prod_{r \in S_0} (\varepsilon^r- \varepsilon^{\zeta(r)}) \prod_{1 \leq l \leq p_0} \varepsilon^{\eta(l)} \mathcal F\right).
\end{equation}

A slight variation on the argument from Theorem \ref{T:|S|geq3} implies that
\begin{multline}
\nu\left(\prod_{\stackrel{S \subseteq [n]}{0 \leq p \leq m}} T_{S,p}^{k\left(S,p\right)}\right) =\frac{1}{N}\left(k(S_0,p_0)-1\right)\varrho_{2p_0} \nu \left(\frac{\mathcal F}{T_{S_0,p_0}}\right)\\+ \frac{1}{N}\sum_{(S_1,p_1) \neq (S_0,p_0)} k(S_1,p_1)\delta_{S_1, S_0}\varrho_{p_1+p_0} \frac{1}{N}\nu \left(\frac{\mathcal F}{T_{S_1,p_1}}\right)\\
+ \beta^2\varrho_{p_0} \nu\left(\mathcal F \left(R_{r_1, r_2}-R_{r_1, \zeta\left(r_2\right)} - R_{\zeta\left(r_1\right), r_2}+R_{\zeta\left(r_1\right), \zeta\left(r_2\right)}\right)\right)
+ O(k+1)
\end{multline} 
where we interpret a term in which a coefficient evaluates to $0$ as being $0$.

To check the induction step, notice that 
\begin{equation}
R_{r_1, r_2}-R_{r_1, \zeta\left(r_2\right)} - R_{\zeta\left(r_1\right), r_2}+R_{\zeta\left(r_1\right), \zeta\left(r_2\right)} = T_{\{r_1, r_2\}}-T_{\{r_1, \zeta\left(r_2\right)\}} - T_{\{\zeta\left(r_1\right),r_2\}}+T_{\{\zeta\left(r_1\right), \zeta\left(r_2\right)\}}.
\end{equation}  
By Lemma \ref{L:Ind} for each pair of replica indices such that $\{\ell, \ell'\} \varDelta \{r_1, r_2\} \neq \varnothing$ 
\begin{equation}
\left|\nu\left(\mathcal F T_{\{\ell, \ell'\}} \right)\right| = O(k+1).
\end{equation}
Thus we have
\begin{multline}
\label{L:Key-Step}
\nu\left(\prod_{\stackrel{S \subseteq [n]}{0 \leq p \leq m}} T_{S,p}^{k\left(S,p\right)}\right) =\frac{1}{N}\left(k(S_0,p_0)-1\right)\varrho_{2p_0} \nu \left(\frac{\mathcal F}{T_{S_0,p_0}}\right)\\+ \frac{1}{N}\sum_{(S_1,p_1) \neq (S_0,p_0)} k(S_1,p_1)\delta_{S_1, S_0}\varrho_{p_1+p_0} \nu \left(\frac{\mathcal F}{T_{S_1,p_1}}\right)\\
+ \beta^2\varrho_{p_0} \nu\left(\mathcal F T_{\{i_1, i_2\}}\right)+ O(k+1).
\end{multline} 

By our assumption on the covariance structure of $\{Y_{S,p}\}$ and the strong induction hypothesis,
\begin{equation}
\nu\left(\prod_{\stackrel{S \subseteq [n]}{0 \leq p \leq m}} T_{S,p}^{k\left(S,p\right)}\right) = N^{-\frac{k}{2}}\mathcal M(S_0) \E\left[\prod_{\stackrel{S \neq S_0}{0 \leq p \leq m}} Y_{S,p}^{k(S,p)}\right]  + O(k+1)
\end{equation} 
where
\begin{multline}
\mathcal M(S_0)= \left(k(S_0,p_0)-1\right)\varrho_{2p_0}\E\left[Y_{S_0,p_0}^{k(S_0, p_0) -2}\prod_{p < p_0} Y_{(S_0, p)}^{k(S_0, p)}\right]\\
+ \sum_{p_1 < p_0} k(S_0,p_1)\varrho_{p_1+p_0} \E\left[Y_{S_0,p_0}^{k(S_0, p_0) -1}Y_{S_0,p_1}^{k(S_0, p_1) -1}\prod_{p \notin \{p_0, p_1\}} Y_{(S_0, p)}^{k(S_0, p)}\right]\\
+ \beta^2\varrho_{p_0}\E\left[Y_{S_0,p_0}^{k(S_0, p_0) -1}\prod_{p < p_0} Y_{(S_0, p)}^{k(S_0, p)}Y_{S_0,0}\right].
\end{multline}
Using the Gaussian structure of the family $\{Y_{S_0,p}\}_{0 \leq p \leq m}$
\begin{multline}
\E\left[Y_{S_0,p_0}^{k(S_0, p_0) -1}\prod_{p < p_0} Y_{(S_0, p)}^{k(S_0, p)}Y_{S_0,0}\right] =\\
 \left(k(S_0,p_0)-1\right)\E\left[Y_{S_0,p_0}Y_{S_0,0}\right]\E\left[\frac{\partial}{\partial Y_{S_0,p_0}} \left(Y_{S_0,p_0}^{k(S_0, p_0) -1}\prod_{p < p_0} Y_{(S_0, p)}^{k(S_0, p)}\right)\right]\\
+\sum_{p_1 \neq p_0} \E\left[Y_{S_0,p_1}Y_{S_0,0}\right]\E\left[\frac{\partial}{\partial Y_{S_0,p_1}} \left(Y_{S_0,p_0}^{k(S_0, p_0) -1}\prod_{p < p_0} Y_{(S_0, p)}^{k(S_0, p)}\right)\right].
\end{multline}
The induction step now follows easily from the covariance structure of $\{Y_{S,p}\}$.
\end{proof}

Before going one step further in our sequence of CLT's we need a preliminary calculation.
\begin{lemma}
\label{L:Ind-Step}
Let $m,n \in \N$ be fixed.
Consider $\mathcal G$ be a monomial of degree $k$ in the truncated multi-overlaps $\{T_{r,p}\}_{r \in [n],0 \leq  p \leq m} \cup \{T_{\varnothing , p}\}_{0 \leq p \leq m}$.  If $S_0= \{r_0, n+1\}$ for some $r_0 \in [n]$ then
\begin{equation}
\nu\left(T_{S_0,p_0}T_{S_0}\mathcal G\right)-\nu\left(T_{S_0,p_0}T_{S_0}\right) \nu\left(\mathcal G\right) = O\left(k+3\right)
\end{equation}
for any $p_0 \in \N$.
\end{lemma}
\begin{proof}
A calculation along the lines of \eqref{L:Key-Step}, bounding errors using Lemma \ref{L:Mom-Bound} gives
\begin{multline}
 \nu\left(T_{S_0,p_0}T_{S_0}\mathcal G\right)= \\
 \frac{1}{N} \varrho_{p_0}\nu(\mathcal G) + \beta^2\varrho_{p_0}\nu\left(T_{S_0}\mathcal G \left(T_{r_0, n+1}-T_{r_0, \zeta\left(n+1\right)} - T_{\zeta\left(r_0\right), n+1}+T_{\zeta\left(r_0\right), \zeta\left(n+1\right)}\right)\right) + O\left(k+3\right). 
\end{multline}
By Lemma \ref{L:Ind} this gives the approximate identity 
\begin{equation}
 \nu\left(T_{S_0,p_0}T_{S_0}\mathcal G\right)= \frac{1}{N} \varrho_{p_0}\nu\left(\mathcal G\right) + \beta^2\varrho_{p_0}\nu\left(T^2_{S_0} \mathcal G\right) + O\left(k+3\right). 
\end{equation}

In particular we may set $p_0=0$ to obtain
\begin{equation}
 \nu\left(T^2_{S_0}\mathcal G\right)= \frac{1}{N}\frac{ \varrho_{0}}{1 - \beta^2 \varrho_0} \nu\left(\mathcal G\right) + O\left(k+3\right). 
\end{equation}
The lemma follows easily.
\end{proof}
\begin{theorem}[Joint Moments: $|S|\geq 1$]
\label{T:|S|geq1}
Consider the family $\{T_{S,p}\}_{S\subseteq [n],\: 0 \leq p \leq m}$ and an arbitrary collection of integers $\{k\left(S,p\right)\}_{S\subseteq [n], \: 0 \leq p \leq m}$.  Suppose that $k(\varnothing, p)= 0$ for all $0 \leq p \leq m$.  Let $k= \sum_{S \subseteq[n], 0 \leq p \leq m} k\left(S,p\right)$.  Then we have
\begin{equation}
\left|\nu\left(\prod_{\stackrel{S \subseteq [n]}{0 \leq p \leq m}} T_{S,p}^{k\left(S,p\right)}\right) - N^{-\ffrac{k}{2}} \prod_{S \subseteq [n]} \E\left[ \prod_{0 \leq p \leq m} Y_{S,p}^{k\left(S,p\right)}\right]\right| \leq \frac{C}{N^{\frac{k+1}{2}}} 
\end{equation}
\end{theorem}
\begin{proof}
Let
\begin{equation}
k_1=\sum_{\stackrel{S: |S| \geq 3}{0 \leq p \leq m}} k(S,p), \hspace{12pt}  k_2= \sum_{\stackrel{|S|=2}{0 \leq p \leq m}} k(S,p)), \hspace{12pt}  k_3=\sum_{\stackrel{|S|=1}{0 \leq p \leq m}} k(S, p)
\end{equation}
and 
\begin{equation}
\tilde k_3=\sum_{\stackrel{|S|=1}{p \neq 1}} k(S, p)
\end{equation}
The proof proceeds in three induction steps.   At the first step we set $k_1=k_2=0$ and proceed by induction on $\tilde k_3$.  This will be detailed below.  Assuming for the moment that we have proved the statement when $k_1, k_2 =0$ and $\tilde k_3$ is arbitrary, an argument analogous  to Theorem \ref{T:|S|geq2} shows that the result holds for $k_1, k_2$ arbitrary. Let us note that there is a bit of work to be done in that the initial step for $k_2$ is not covered by Theorem \ref{T:CLT-Tal}.  This can be taken care of using techniques along the lines of Lemma \ref{L:Ind-Step}.

Let us consider the first induction, in which $k_1=k_2=0$.
We identify moments by induction on $\tilde k_3$.  The induction hypothesis we use is similar to that of Theorem \ref{T:|S|geq2}.

Consider a collection of exponents $\{k(S, p)\}_{|S|=1,\: 0 \leq p \leq m}$.
Let $(S_0, p_0)$ be chosen so that $p_0$ is maximal and $k(S_0, p_0) \neq 0$.  If there are multiple choices of $S_0$ we choose arbitrarily among those with $k(S_0,p_0)$ maximal as well.  As an induction hypothesis for $\tilde k_3$, let us suppose that we have the statement for all choices of exponents $\{k^*(S, p)\}_{|S| =  1, 0 \leq  p \leq m}$ with $\{k^*(S,1)\}_{|S| = 1}$ arbitrary and  so that
\begin{equation}
k^*(S,p) \leq k(S,p),  \hspace{12pt} k^*(S_0,p_0) < k(S_0,p_0)
\end{equation}
for $p \neq 1$.

The initial step, when $\tilde k_3=0$ and the remaining exponents are arbitrary, is part of the content of Theorem \ref{T:CLT-Tal}.  We thus proceed to the induction step.
Let us denote $S_0=\{r_0\}$ and when we want to be explicit
\begin{equation}
 T_{S_0,p}= T_{r_0,p}.
\end{equation}
For convenience let
\begin{equation}
\mathcal{F}= T_{r_0,p_0}^{k(S_0, p_0)-1} \prod_{(S,p) \neq (r_0,p_0)} T_{S,p}^{k(S,p)}.
\end{equation}
By symmetry, we have
\begin{equation}
\nu\left(\prod_{\stackrel{S \subseteq [n]}{0 \leq p \leq m}} T_{S,p}^{k\left(S,p\right)}\right) = \nu\left((\varepsilon^{r_0}- \varepsilon^{\zeta(r_0)}) \prod_{l=1}^{p_0} \varepsilon^{\eta(l)} \mathcal F\right).
\end{equation}

A slight variation on the argument from Theorem \ref{T:|S|geq2}  using the symmetry of replicas implies that
\begin{multline}
\nu\left(\prod_{\stackrel{S \subseteq [n]}{0 \leq p \leq m}} T_{S,p}^{k\left(S,p\right)}\right) =\frac{1}{N}\left(k(r_0,p_0)-1\right)\pi_{2p_0-1} \nu \left(\frac{\mathcal F}{T_{r_0,p_0}}\right) \\
+ \frac{1}{N}\sum_{p_1 < p_0} k(r_0,p_1)\pi_{p_1+p_0-1} \frac{1}{N}\nu \left(\frac{\mathcal F}{T_{r_0,p_1}}\right)\\
+ \beta^2 \pi_{p-2}\sum_{1 \leq l \leq p} \nu\left(\mathcal F\left(R_{r_0, \eta \left(l\right)}-R_{\zeta\left(r_0\right), \eta\left(l\right)}\right)\right)\\
\hspace{48pt} + \beta^2 \pi_{p}\sum_{ \ell \notin\{r_0,\zeta\left(r_0\right)\}\cup \eta \left([p]\right)} \nu\left(\frac{\mathcal F}{T_{S,p}(\ell)}T_{S,p}(\ell) \left(R_{r_0, \ell}-R_{\zeta\left(r_0\right), \ell}\right)\right)\\
 \hspace{48pt} - \beta^2 d \pi_{p}\nu\left(\mathcal F\left(R_{r_0, z}-R_{\zeta\left(r_0\right), z}\right)\right)
+ O(k+1)
\end{multline} 
where $z$ denotes a replica disjoint from those introduced thus far and 
\begin{equation}
d = \left|\{r \in [n] : \text{ $k(\{r\},p) \neq 0$ for some $p$}\}\right| + \sum_{(S,p)} k(S,p)(p+1)
\end{equation}
and we interpret a term in which a coefficient evaluates to $0$ as being $0$.

Now, we have the identities used previously in Lemma \ref{L:|S|=1}:
\begin{align*}
T_{r,p}(\ell) = &T_{r, p} - T_{\ell, p}\\
R_{r_0,\ell}-R_{\zeta(r_0), \ell} = & T_{\{r_0, \ell\}}-T_{\{\zeta(r_0), \ell\}} + T_{r_0, 1}-T_{\zeta(r_0), 1}
\end{align*}
if $\ell$ is the image of an $\alpha$
and 
\begin{align*}
T_{r,p}(\ell) = &T_{\{r, \ell\}, p-1} + T_{r, p}\\
R_{r_0,\ell}-R_{\zeta(r_0), \ell} = & T_{\{r_0, \ell\}}-T_{\{\zeta(r_0), \ell\}} + T_{r_0, 1}-T_{\zeta(r_0), 1}
\end{align*}
if $\ell$ is in the image of a $\gamma$.

Noting that the overlap involves a free index in the remaining summands, inserting these two identities into their respective summands and applying Lemma \ref{L:Ind} whenever possible gives
\begin{align*}
&\nu\left(\mathcal F\left(R_{r_0, \eta\left(l\right)}-R_{\zeta\left(r_0\right), \eta\left(l\right)}\right)\right)= \nu\left(\mathcal F\left(R_{r_0, z}-R_{\zeta\left(r_0\right), z}\right)\right)=\nu\left(\mathcal F T_{r_0, 1}\right) + O(k+1),\\
&\nu\left(\frac{\mathcal F}{T_{S,p}(\ell)}T_{S,p}(\ell)\left(R_{r_0, \ell}-R_{\zeta\left(r_0\right), \ell}\right)\right) =\nu\left(\mathcal F T_{r_0, 1}\right) + O(k+1) \text{ if $\ell$ is the image of an $\alpha$,}\\
&\nu\left(\frac{\mathcal F}{T_{S,p}(\ell)}T_{S,p}(\ell)\left(R_{r_0, \ell}-R_{\zeta\left(r_0\right), \ell}\right)\right)= \nu\left(\mathcal F T_{r_0,1}\right) +\delta_{S, \{r_0\}}\nu\left(\frac{\mathcal F}{T_{S,p}} T_{\{S, \ell\},p-1}T_{\{r_0,\ell\}}\right)\\& \hspace{200pt}+ O(k+1)\text{ if $\ell$ is the image of a $\gamma$}.
\end{align*}
Collecting terms we have
\begin{multline}
\nu\left(\prod_{\stackrel{S \subseteq [n]}{0\leq p \leq m}} T_{S,p}^{k\left(S,p\right)}\right) =\\
\frac{1}{N}\left(k(r_0,p_0)-1\right)\pi_{2p_0-1} \nu \left(\frac{\mathcal F}{T_{r_0,p_0}}\right) + \frac{1}{N}\sum_{p_1 < p_0} k(r_0,p_1)\pi_{p_1+p_0-1} \frac{1}{N}\nu \left(\frac{\mathcal F}{T_{r_0,p_1}}\right) \\
+ \beta^2\left\{p_0 \pi_{p_0-2} - (p_0+2)\pi_{p_0}\right\}\nu\left(\mathcal F T_{r_0, 1}\right) + \beta^2 \pi_{p_0}(k(r_0, p_0) -1)p_0 \nu\left(\frac{\mathcal F}{T_{r_0,p_0}}T_{\{r_0, z\},p_0-1}T_{\{r_0, z\}, 1}\right)\\
\hspace{48pt} + \beta^2 \pi_{p_0}\sum_{p_1<p_0  } k(r_0, p_1)p_1\nu\left(\frac{\mathcal F}{T_{r_0,p_1}}T_{\{r_0, z\},p_1-1}T_{\{r_0, z\}, 1}\right) + O(k+1)
\end{multline} 
where as usual we interpret terms with $0$ coefficient as being $0$ and $z$ is a new index.

The induction hypothesis and Lemma \ref{L:Ind-Step} imply
\begin{align*}
\nu\left(\mathcal F T_{r_0, 1}\right) = & \sum_{p_1 \leq p_0} \nu\left(T_{r_0, p_1}T_{r_0,1}\right) \nu\left(\frac{\partial}{\partial T_{r_0, p_1}} \mathcal F\right)+ O(k+1).\\
k(r_0, p_1)\nu\left(\frac{\mathcal F}{T_{r_0,p_1}}T_{\{r_0, z\},p_1-1}T_{\{r_0,z\}}\right) =&\nu\left(T_{\{r_0, z\},p_1-1}T_{\{r_0,z\}}\right) \nu\left(\frac{\partial}{\partial T_{r_0, p_1}} \mathcal F\right)  + O(k+1).
\end{align*}
Gathering these observations together,
\begin{multline}
\nu\left(\prod_{\stackrel{S \subseteq [n]}{0 \leq p \leq m}} T_{S,p}^{k\left(S,p\right)}\right) =B(p_0,p_0)\nu\left(\frac{\partial}{\partial T_{r_0, p_0}} \mathcal F\right) + \sum_{p_1 < p_0} B(p_0,p_1) \nu\left(\frac{\partial}{\partial T_{r_0, p_1}} \mathcal F\right) + O(k+1).
\end{multline} 
where 
\begin{equation}
B(p, \tilde p)=\frac{1}{N}\pi_{\tilde p+ p-1}+  \beta^2\left[p \pi_{p-2} - (p+2)\pi_{p}\right]\nu\left(T_{r_0, \tilde p}T_{r_0, 1}\right)\\
+\beta^2\tilde p \pi_{p} \nu\left(T_{\{r_0, z\},\tilde p-1}T_{\{r_0, z\}, 1}\right)
\end{equation}
The result now follows from Lemma \ref{L:|S|=1} and the induction hypothesis.
\end{proof}
\begin{lemma}
\label{L:Ind-Step-1}
Let $m \in \N$ be fixed and $\mathcal G$ be a monomial of degree $k$ in the truncated multi-overlaps $\{T_{\varnothing , p}\}_{0 \leq p \leq m}$. 
Then
\begin{equation}
\nu\left(T_{1, p_0}T_{1,1}\mathcal G\right)-\nu\left(T_{1, p_0}T_{1, 1}\right) \nu\left(\mathcal G\right) = O\left(k+3\right)
\end{equation}
fo each $p_0 \in \N$.
\end{lemma}
\begin{proof}
The statement follows from slight modifications to the argument of Lemma \ref{L:Ind-Step}.
\end{proof}
\begin{theorem}[Joint Moments: The Full Picture]
Consider the family $\{T_{S,p}\}_{S\subseteq [n],\:0 \leq  p \leq m}$ and an arbitrary collection of integers $\{k\left(S,p\right)\}_{S\subseteq [n], \:0 \leq  p \leq m}$.  Let $k= \sum_{(S,p)} k\left(S,p\right)$.  Then we have
\begin{equation}
\left|\nu\left(\prod_{(S,p)} T_{S,p}^{k\left(S,p\right)}\right) - N^{-\ffrac{k}{2}} \prod_{S \subseteq [n]} \E \left[\prod_{p=0}^m Y_{S,p}^{k\left(S,p\right)}\right]\right| \leq \frac{C}{N^{\frac{k+1}{2}}} 
\end{equation}
\end{theorem}
\begin{proof}
Let us define
\begin{equation}
k_1 = \sum_{(S,p): \: |S| \geq 1} k(S,p), \quad k_2 = \sum_{0 \leq p \leq m}k(\varnothing, p).
\end{equation}
Similar to our previous arguments, we may reduce to the case $k_1=0$.  Now the canonical pair is $(\varnothing , 2)$ and we will prove the remaining identifications by induction on
\begin{equation}
\tilde k_2 = \sum_{p \neq 2 }k(\varnothing, p).
\end{equation}

To state our induction hypothesis, let $p_0$ be maximal so that $p_0 \neq 2$ and $k(\varnothing, p_0) \neq 0$.  Suppose that we have the statement for all choices of exponents $\{k^*(\varnothing, p)\}_{0 \leq  p \leq m}$ with $\{k^*(\varnothing,2)\}$ arbitrary and  so that
\begin{equation}
k^*(\varnothing,p) \leq k(\varnothing,p),  \hspace{12pt} k^*(\varnothing, p_0) < k(\varnothing,p_0)
\end{equation}
for $p \neq 2$.
The case $\tilde k_2 = 0$ and $k(\varnothing, 2)$ arbitrary is part of Theorem \ref{T:CLT-Tal}.  
 
To check the induction step at $\tilde k_2 =K>0$, let us fix a monomial 
\begin{equation}
\prod_{1 \leq p \leq m} T_{S, p}^{k(\varnothing , p)}
\end{equation}
with $\tilde k_2$ degree $K$ and so that $p_0$ is maximal with $p_0 \neq 2$ and $k(\varnothing , p_0) \neq 0$.  Let us denote the total degree of this monomial by $k$.
For notational convenience let us define the monomial $\mathcal F$ by
\begin{equation}
\prod_{1 \leq p \leq m} T_{S, p}^{k(\varnothing , p)}= T_{\varnothing, p_0}\mathcal F.
\end{equation}

As usual we have by symmetry, 
\begin{multline}
\nu\left(T_{\varnothing, p_0}\mathcal F\right) = \frac{1}{N}\left(k(\varnothing,p_0)-1\right) \left(q_{2p_0} - 2 q_{p_0}\right) \nu \left(\frac{\mathcal F}{T_{\varnothing, p_0}}\right) \\
+ \frac{1}{N}\sum_{p \neq p_0} k(\varnothing ,p) \left(q_{p_0+p} - q_{p_0}q_p\right) \frac{1}{N}\nu \left(\frac{\mathcal F}{T_{\varnothing, p}}\right) + \nu\left(\left(\prod_{l=1}^{p_0} \varepsilon^{\eta(l)} - q_{p_0}\right)\mathcal F^-\right).
\end{multline}

Now we may expand this product using a total of $d=\sum_{1 \leq p \leq m} k(\varnothing, p) p$ replicas.
By Lemmas \ref{L:cavitation}, \ref{L:Cav-Bound} and \ref{L:Mom-Bound}
\begin{multline}
\nu\left(\left(\prod_{l=1}^{p_0} \varepsilon^{\eta\left(l\right)}- q_{p_0}\right)\mathcal F^{-}\right)=
\beta^2 \left(q_{{p_0}-2}- q_{p_0}q_2\right){p_0 \choose 2} \nu\left(\mathcal FT_{\varnothing, 2}\right)\\
+  \beta^2 \frac{d\left(d+1\right)}{2} \left(q_{{p_0}+2}-q_{p_0} q_2\right)\nu\left(\mathcal FT_{\varnothing, 2}\right)\\
- \beta^2 p_0 d \left(q_{p_0}-q_{p_0}q_2\right) \nu\left(\mathcal F T_{\varnothing, 2}\right)\\
- \beta^2 d \left(q_{{p_0}+2}-q_{p_0}q_2\right)\sum_{\ell \notin \eta([p_0])} \nu\left(\mathcal F\left(T_{\ell, 1} + T_{\varnothing, 2}\right)\right)\\
+\beta^2 \left(q_{p_0}- q_{p_0} q_2\right)\sum_{\stackrel{\ell \in \eta([p_0])}{ \ell' \notin \eta([p_0])}} \nu\left(\mathcal F\left(T_{\{\ell, \ell'\}}+ T_{\ell', 1} + T_{\varnothing, 2}\right)\right)\\
+\beta^2 \left(q_{p_0+2}- q_{p_0} q_2\right)\sum_{\{\ell, \ell'\} \cap \eta([p_0]) = \varnothing} \nu\left(\mathcal F\left(T_{\{\ell, \ell'\}}+ T_{\ell, 1}+ T_{\ell', 1} + T_{\varnothing, 2}\right)\right) + O(k+1)
\end{multline}

Consider terms four and five.   We write $\mathcal F= \frac{\mathcal F}{T_{\varnothing, p}(\ell)} T_{\varnothing, p}(\ell)$.
Since 
\begin{equation}
T_{\varnothing, p}(\ell) = T_{\ell, p-1} + T_{\varnothing, p}
\end{equation}
we have
\begin{align*}
\nu\left(\mathcal F\left(T_{\ell, 1}+ T_{\varnothing, 2}\right)\right)= & \nu\left(T_{1, p-1}T_{1,1}\right)  \nu\left(\frac{\mathcal F}{T_{\varnothing, p}} \right) + \nu\left(T_{\varnothing, 2}\mathcal F \right) + O(k+1)\\
\nu\left(\mathcal F\left(T_{\{\ell, \ell'\}}+ T_{\ell', 1} + T_{\varnothing, 2}\right)\right)= & \nu\left(T_{1, p-1}T_{1,1}\right)  \nu\left(\frac{\mathcal F}{T_{\varnothing, p}} \right) + \nu\left(T_{\varnothing, 2}\mathcal F \right)  + O(k+1)
\end{align*}
where we have used Lemmas \ref{L:Ind}, \ref{L:Ind-Step} and \ref{L:Ind-Step-1}.

For the sixth term, there are two cases to distinguish: either the replicas $\ell, \ell'$ appear in the same truncated overlap or in distinct truncated overlaps.  By Lemmas \ref{L:Ind}, \ref{L:Ind-Step} and  \ref{L:Ind-Step-1} 
\begin{multline}
\nu\left(\mathcal F\left(T_{\{\ell, \ell'\}} + T_{\ell,1} + T_{\ell' , 1} + T_{\varnothing, 2}\right)\right) = 
\left[\nu\left(T_{\{1, 2\}, p-2}T_{\{1, 2\}}\right) +2\nu\left(T_{1, p-1}T_{1,1}\right) \right]\nu\left(\frac{\mathcal F}{T_{\varnothing, p}}\right) \\
+\nu\left(\mathcal F T_{\varnothing,2}\right) +O(k+1) 
\end{multline}
if $\ell, \ell'$ come from the same overlap and
\begin{multline}
\nu\left(\mathcal F\left(T_{\{\ell, \ell'\}} + T_{\ell,1} + T_{\ell' , 1} + T_{\varnothing, 2}\right)\right)= \nu\left(T_{1, p-1}T_{1,1}\right)\nu\left(\frac{\mathcal F}{T_{\varnothing, p}}\right) \\
+ \nu\left(T_{1, p'-1}T_{1,1}\right)\nu\left(\frac{\mathcal F}{T_{\varnothing, p'}}\right)+ \nu\left(\mathcal FT_{\varnothing, 2}\right)  + O(k+1)
\end{multline}
if $\ell, \ell'$ come from distinct overlaps $T_{\varnothing, p}(\ell), T_{\varnothing, p'}(\ell')$ respectively.

All that remains is to collect terms.  It is most convenient to collect the contributions coming from a fixed replica first.  For each $\ell \notin \eta([p_0])$, the term $\nu\left(\frac{\mathcal F}{T_{\varnothing, p}(\ell)}\right)$ appears with the coefficient
\begin{multline}
\beta^2 \left\{  \left(q_{p_0+2}- q_{p_0} q_2\right)\left(-1 -p_0\right) + \left(q_{p_0}- q_{p_0} q_2\right) p_0\right\}  \nu\left(T_{1, p-1} T_{1,1}\right) \\
+ \frac{p-1}{2} \left(q_{p_0+2}- q_{p_0} q_2\right)\nu\left(T_{\{1, 2\}, p-2}T_{\{1, 2\}}\right).
\end{multline}
where we have counted the term $\nu\left(T_{\{1, 2\}, p-2}T_{\{1, 2\}}\right)$ contributed by the pair $\{\ell, \ell'\}$ as a one half contribution to each of $\nu\left(\frac{\mathcal F}{T_{\varnothing, p}(\ell)}\right)$, $\nu\left(\frac{\mathcal F}{T_{\varnothing, p}(\ell')}\right)$.

On the other hand,  $\nu\left(\mathcal F T_{\varnothing, 2}\right)$ appears with coefficient
\begin{multline}
\beta^2 \lambda_{p_0}=\beta^2 {p_0 \choose 2}\left(q_{{p_0}-2}- q_{p_0}q_2\right)
+  \beta^2 {p_0+1 \choose 2}\left(q_{{p_0}+2}-q_{p_0} q_2\right)
- \beta^2 p_0^2 \left(q_{p_0}- q_{p_0} q_2\right).
\end{multline}
Using symmetry of replicas in the $\nu$ measure we thus have
\begin{multline}
\nu\left(T_{\varnothing, p_0}\mathcal F\right) = \left(k(\varnothing,p_0)-1\right)\mathcal M^*(p_0) \nu \left(\frac{\mathcal F}{T_{\varnothing, p_0}}\right)\\
+ \sum_{p \neq p_0} k(\varnothing, p)\mathcal M(p)\nu \left(\frac{\mathcal F}{T_{\varnothing, p}}\right) + \beta^2 \lambda_{p_0}\nu\left(T_{\varnothing, 2} \mathcal F\right) + O(k+1).
\end{multline}
where
\begin{equation}
\mathcal M^*(p_0)=\\ \frac{1}{N} \left(q_{2p_0} - 2 q_{p_0}\right) + \kappa_{p_0, p_0} \nu\left(T_{1, p_0-1}T_{1,1}\right) +\omega_{p_0, p_0}  \nu\left(T_{\{1, 2\} p_0-2}T_{\{1, 2\}}\right) 
\end{equation}
and
\begin{equation}
\mathcal M(p)= \frac{1}{N} \left(q_{p+p_0} - q_{p}q_{p_0}\right) + \kappa_{p_0, p} \nu\left(T_{1, p-1}T_{1,1}\right) +\omega_{p_0, p}  \nu\left(T_{\{1, 2\} p-2}T_{\{1, 2\}}\right)
\end{equation}
where, as usual, if a coefficient is $0$ then we interpret the summand as being $0$.

By the (strong) induction hypothesis,
\begin{equation}
 \nu\left(T_{\varnothing, 2} \mathcal F\right)= \sum_{1 \leq p \leq m} \nu\left(T_{\varnothing, 2} T_{\varnothing, p}\right)\nu\left(\frac{\partial}{\partial T_{\varnothing, p}} \mathcal F\right)+ O(k+1).
\end{equation}
Of course $k(\varnothing , p) \nu\left(\frac{\mathcal F}{T_{\varnothing, p}}\right) = \nu\left(\frac{\partial}{\partial T_{\varnothing, p}} \mathcal F \right)$, so collecting terms, comparing them to Lemma \ref{L:|S|=0} and applying the induction hypothesis to get
\begin{equation}
\nu\left(\frac{\partial}{\partial T_{\varnothing, p'}} \mathcal F \right)= \E\left[\frac{\partial}{\partial Y_{\varnothing, p'}} \left(Y_{\varnothing p_0}^{k(\varnothing, p_0)-1} \prod_{p \neq p_0} Y_{\varnothing, p}^{k(\varnothing, p)}\right) \right] + O(k-1)
\end{equation}
finishes the induction step.

\end{proof}

\section[Appendix]

Here we provide proof of the main tool for estimation of the errors incurred in the application of Lemma \ref{L:Cav-Bound}.
We shall use the inequality for overlaps due to Talagrand mentioned in the introduction.

\begin{lemma}
\label{L:Mom-Bound-2}[Theorem 2.5.1 \cite{Talagrand-book}]
There exists a $\beta_0 > 0$ such that for all $\beta \leq \beta_0$,
\begin{equation}
\nu\left(\left(R_{1,2}- q_2\right)^{2k}\right) \leq \left(\frac{Lk}{N}\right)^k
\end{equation}
for all $k \in \N$.  Here $L>0$ denotes a constant which is independent of $k, N$.
\end{lemma}

\begin{lemma}
\label{L:Mom-Bound-1}
Let $\beta_0$ be as in Lemma \ref{L:Mom-Bound-2}.  Suppose $S \subset \mathbb N$ is finite and let $p \in\N$ be fixed.  Consider the truncated overlap $T_{S,p}$.  For all $\beta \leq \beta_0$ and any $k\in \N$, there exists a constant $C$ depending on $S, p , k$ ( and $\beta$) so that
\begin{equation}
\nu \left(T_{S,p}^{2k}\right) \leq \frac{C(S, p, k)}{N^{k}}
\end{equation}
and
\begin{equation}
\nu \left(\left(T_{S,p}^-\right)^{2k}\right) \leq \frac{C(S, p, k)}{N^{k}}.
\end{equation}
\end{lemma}

\begin{proof}
We prove these inequalities simultaneously by induction on $k$.  
Let us note that we are not after exponential moments here.  It does not seem clear that we obtain exponential moments uniformly in $\beta_0$ over arbitrary choice of pairs $(S, p)$ without a more careful consideration of the interpolation term.  This would require more delicate analysis than we need, and so will not be pursued here.

Consider the case $k=1$.
Using the symmetry of sites,
\begin{align*}
\nu\left(T^2_{S,p}\right)&= \nu\left(\prod_{r\in S}\left(\varepsilon^r - \varepsilon^{\zeta(r)}\right)\prod_{l=1}^p\varepsilon^{\eta(l)} T_{S,p}\right) \\=&
\frac{1}{N}\nu\left(\prod_{r\in S}\left(\varepsilon^r - \varepsilon^{\zeta(r)}\right) \prod_{l=1}^p\varepsilon^{\eta(l)} \prod_{r\in S}\left(\varepsilon^r - \varepsilon^{\alpha(r)}\right) \prod_{l=1}^p \varepsilon^{\gamma(l)}\right) +  \nu\left(\prod_{r\in S}\left(\varepsilon^r - \varepsilon^{\zeta(r)}\right) T_{S,p}^-\right) 
\end{align*}
On the other hand, Lemma \ref{L:Cav-Bound} implies that
we have
\begin{equation}
\label{Eq:-1}
\left|\nu\left(\prod_{i\in S}\left(\varepsilon^r - \varepsilon^{\zeta(r)}\right) \prod_{l=1}^p\varepsilon^{\eta(l)}  T_{S,p}^- \right) \right| \leq \frac{2^{|S|}}{\sqrt N} \nu\left(\left(T_{S,p}^-\right)^2 \right)^{\frac{1}{2}} \leq  2\frac{2^{|S|}}{\sqrt N} \nu \left(T_{S,p}^2 \right)^{\frac{1}{2}} + 2\frac{4^{|S|}}{N^{\frac{3}{2}}} 
\end{equation}

This implies the bound
\begin{equation}
\nu\left(T^2_{S,p}\right) \leq \frac{4^{|S|}}{N} +  2\frac{2^{|S|}}{\sqrt N} \nu \left(T_{S,p}^2 \right)^{\frac{1}{2}} + 2\frac{4^{|S|}}{N^{\frac{3}{2}}} 
\end{equation}
It follows easily that
\begin{equation}
\nu\left(T^2_{S,p}\right) \leq \frac{4^{|S|+2}}{N}
\end{equation}
for all $N$ sufficiently large.  The bound for $T_{S,p}^-$ follows from \eqref{Eq:-1}.

For the induction step suppose that we have proved the statement for all $k \leq m$.  Proceeding as above,
\begin{equation}
\nu\left(T^{2m+2}_{S,p}\right) = \nu\left(\prod_{r\in S}\left(\varepsilon^r - \varepsilon^{\zeta(r)}\right) \prod_{l=1}^p\varepsilon^{\eta(l)} T^{2m+1}_{S,p}\right).
\end{equation}
Now we may write
\begin{equation}
\label{Exp1}
T^{2m+1}_{S,p} = \left(T_{S,p}^-\right)^{2m+1} +\mathcal E
\end{equation}
where 
\begin{equation}
\left|\mathcal E\right| \leq \sum_{s=0}^{2m} {2m+1 \choose s} \left(\left|T_{S,p}^-\right|\right)^s\frac{2^{|S|(2m+1-s)}}{N^{2m+1-s}}.
\end{equation}
By H\"{o}lder's inequality and the induction hypothesis
\begin{equation}
\nu\left(\left(\left|T_{S,p}^-\right|\right)^s\right)\leq\nu\left(\left(T_{S,p}^-\right)^{2m}\right)^{\frac{s}{2m}} \leq \frac{C}{N^{\ffrac{s}{2}}}
\end{equation}
so that 
\begin{equation}
\nu\left(|\mathcal E| \right) \leq \frac{(2m+1)2^{|S|}}{N} \left(\frac{C}{\sqrt N} + \frac{2^{|S|}}{N}\right)^{2m}
\end{equation}
By Proposition 2.4.7 of \cite{Talagrand-book}
\begin{multline}
\left|\nu\left(\prod_{r\in S}\left(\varepsilon^r - \varepsilon^{\zeta(r)}\right) \prod_{l=1}^p\varepsilon^{k(l)} \left(T_{S,p}^-\right)^{2m+1}\right) - \nu_0\left(\prod_{r\in S}\left(\varepsilon^r - \varepsilon^{\zeta(r)}\right) \prod_{l=1}^p\varepsilon^{\eta(l)}  \left(T_{S,p}^-\right)^{2m+1}\right)\right| \\\leq 2 d^2 \beta^2 \exp\left(4 d^2 \beta^2\right)\nu \left(\left(T_{S, p}^-\right)^{2m+2}\right)^{\frac{2m+1}{2m+2}} \nu\left(\left(R_{1,2} -q_2\right)^{2m+2}\right)^{\frac{1}{2m+2}}
\end{multline}
where $d=(2m+2)\left(|S|+p\right) + |S| $.

Since 
\begin{equation}
\nu_0\left(\prod_{r\in S}\left(\varepsilon^r - \varepsilon^{\zeta(r)}\right) \prod_{l=1}^p\varepsilon^{\eta(l)}  \left(T_{S,p}^-\right)^{2m+1}\right) = 0
\end{equation}
we have 
\begin{multline}
\nu\left(T_{S,p}^{2m+2}\right) \leq  2 d^2 \beta^2 \exp\left(4 d^2 \beta^2\right)\nu \left(\left(T_{S, p}^-\right)^{2m+2}\right)^{\frac{2m+1}{2m+2}} \nu\left(\left(R_{1,2} -q_2\right)^{2m+2}\right)^{\frac{1}{2m+2}} \\
+ \left(\frac{C m}{N} \right)^{m+1} 
\end{multline}
To get a workable inequality we make one more observation.  Similar to \eqref{Exp1} we have
\begin{equation}
\label{relate}
\left(T_{S,p}^-\right)^{2m+2} = T_{S,p}^{2m+2} + \tilde {\mathcal E}
\end{equation}
where
\begin{equation}
\left|\tilde {\mathcal E} \right | \leq \sum_{s=0}^{2m+1} {2m+2 \choose s} \left(\left|T_{S,p}\right|\right)^s\frac{2^{|S|(2m+2-s)}}{N^{2m+2-s}}.
\end{equation}
Applying the knowledge (from Lemma \ref{L:Mom-Bound-2}) that
\begin{equation}
\nu\left(\left(R_{1,2}-q_2\right)^{2m+2}\right)^\frac{1}{2m+2} \leq \frac{L}{\sqrt N}
\end{equation}
and Minkowski's inequality gives
\begin{equation}
\nu\left(T_{S,p}^{2m+2}\right) \leq  2 d^2 \beta^2 \exp\left(4 d^2 \beta^2\right)\left[ \nu \left(T_{S, p}^{2m+2}\right)^{\ffrac{1}{2m+2}}+ \frac{C'}{N}\right]^{2m+1} \frac{L}{\sqrt N} 
+ \left(\frac{C m}{N} \right)^{m+1} 
\end{equation}
where $C'$ is some appropriately chosen constant.  The first statement of the induction now follows.  To see the second statement, we refer the reader to \eqref{relate}.
\end{proof}

\end{document}